\documentclass{amsart}
\usepackage{a4}
\usepackage{pstricks,pst-node}
\usepackage{xypic}

\numberwithin{equation}{section}
\newcommand{\CRes}{{\bf \sf C}}

\newcommand{\aaf}{\mathfrak a}
\newcommand{\bbf}{\mathfrak b}
\newcommand{\rr}{\mathfrak r}

\newcommand{\N}{\mathbb{N}}

\newcommand{\Bc}{\mathcal{B}}

\newcommand{\Llc}{\mathcal{L}}
\newcommand{\Mmc}{\mathcal{M}}

\newcommand{\SM}{\rm{MinGen}}
\newcommand{\nbc}{{\bf nbc}}
\newcommand{\propP}{{\bf (P)}}
\newcommand{\propH}{{\bf (H)}}
\DeclareMathOperator{\Hilb}{Hilb}  
\DeclareMathOperator{\id}{id}  
\DeclareMathOperator{\ini}{in_{\prec}}

\DeclareMathOperator{\lcm}{lcm}

\DeclareMathOperator{\supp}{supp}  
\DeclareMathOperator{\Tor}{Tor}

\newcommand{\on}[1]{\operatorname{#1}}

\newcommand{\dotcup}{\stackrel{_{_\cdot}}{\cup}}

\newcommand{\poinm}[2]{P^{#2}_{#1}(\ul{x},-t)}
\newcommand{\poinr}[2]{P^{#2}_{#1}(\ul{x},t)}

\newcommand{\hilb}[1]{\Hilb_{#1}(\ul{x},t)}

\DeclareMathOperator{\subg}{\unlhd}

\DeclareMathOperator{\tensor}{\otimes}  
\DeclareMathOperator{\iso}{\cong}  
  
\DeclareMathOperator{\Dirsum}{\bigoplus}  
  
\DeclareMathOperator{\pnt}{\raise 0.5mm \hbox{\large\bf.}}

\theoremstyle{plain}
\newtheorem{Thm}{\bf Theorem}[section]
\newtheorem{Lem}[Thm]{\bf Lemma}  
\newtheorem{Cor}[Thm]{\bf Corollary}  
\newtheorem{Prop}[Thm]{\bf Proposition}  
\newtheorem{Conj}[Thm]{\bf Conjecture} 
\theoremstyle{definition}
  
\newtheorem{Ex}[Thm]{\bf Example}  
  
\newtheorem{Def}[Thm]{\bf Definition}

\newcommand{\D}{\displaystyle}
\newcommand{\sst}{\scriptstyle}

\newcommand{\ul}[1]{\underline{#1}}

\newcommand{\disbreak}{\allowdisplaybreaks}

\newcommand{\lbreak}{$\:$\\}

\begin{document}

\title[Multigraded Hilbert and Poincar\'e-Betti Series]
{On the Multigraded Hilbert and Poincar\'e-Betti Series and the Golod Property of 
Monomial Rings}

\date{\today}

\author[Michael J\"ollenbeck]{Michael J\"ollenbeck$^1$}
\footnotetext[1]{Supported by EC's IHRP program through grant 
HPRN-CT-2001-00272}
\address{Philipps-Universit\"at Marburg \\
Fachbereich Mathematik und Informatik \\
D-35032 Marburg \\ Germany  }
\email{joella@mathematik.uni-marburg.de}



%
%
\begin{abstract}
In this paper we study the multigraded Hilbert and Poincar\'e-Betti series of 
$A=S/\aaf$, where $S$ is the ring of polynomials in $n$ indeterminates 
divided by the monomial ideal $\aaf$. There is a conjecture about the 
multigraded Poincar\'e-Betti series by Charalambous and Reeves which they proved 
in the case, where the Taylor resolution is minimal.
We introduce a conjecture about the minimal $A$-free resolution of the residue
class field and show that this conjecture implies the conjecture of 
Charalambous and Reeves and,in addition, gives a formula for the Hilbert series.
Using Algebraic Discrete Morse theory, we prove that the homology of the 
Koszul complex of $A$ with respect to $x_1,\ldots, x_n$ is isomorphic to a 
graded commutative ring of polynomials over certain sets in the Taylor 
resolution divided by an ideal $\rr$ of relations. This leads to a proof of 
our conjecture for some classes of algebras $A$.
We also give an approach for the proof of our conjecture via Algebraic Discrete
Morse theory in the general case.\\
The conjecture implies that $A$ is Golod if and only if the product 
(i.e. the first Massey operation) on the Koszul homology is trivial.
Under the assumption of the conjecture we finally prove that a very simple
purely combinatorial condition on the minimal monomial generating system 
of $\aaf$ implies Golodness for $A$.
\end{abstract}

%
%
\maketitle

%
%
%
\section{Introduction}\label{sec_intro}
In this note, we study the multigraded 
Hilbert and Poincar\'e-Betti series of algebras $A=S/\aaf$, where $S$ is the 
commutative polynomial ring in $n$ indeterminates and $\aaf$ is a monomial 
ideal with minimal monomial generating system $\SM(\aaf):=\{m_1,\ldots, m_l\}$.

Recall that the multigraded Poincar\'e-Betti series $\poinr{k}{A}$ and $\hilb{A}$ of 
$A$ are defined as
\begin{eqnarray*}
\poinr{k}{A}&:=&\sum_{i=0}^\infty\sum_{\alpha\in\N^n} 
     \dim_k(Tor_i^A(k,k)_\alpha)~\ul{x}^\alpha~t^i,\\
\hilb{A}&:=&\sum_{i=0}^\infty\sum_{\alpha\in\N^n\atop|\alpha|=i} 
     \dim_k(A_\alpha)~\ul{x}^\alpha~t^i.
\end{eqnarray*}

In \cite{CH} Charalambous and Reeves proved that in the case 
where the Taylor resolution of $\aaf$ over $S$ is minimal the Poincar\'e-Betti 
series takes the following form:
\begin{eqnarray*}
  \poinr{k}{A}=\frac{\D\prod_{i=1}^n(1+x_i\;t)}
  {\D 1+\sum_{I\subset\{1,\ldots,l\}}(-1)^{cl(I)}m_I\;
    t^{cl(I)+|I|}},
\end{eqnarray*}
where $cl(I)$ is the number of equivalence classes of $I$ with respect to the 
relation defined as the transitive closure of 
$i\sim j:\Leftrightarrow \gcd(m_i,m_j)\neq 1$ 
and $m_I:=\lcm(m_i\mid i\in I)$ is the least common multiple.\\
In the general case, they conjecture that 
\begin{eqnarray*}
  \poinr{k}{A}=\frac{\D\prod_{i=1}^n(1+x_i\;t)}
  {\D 1+\sum_{I\subset[l] \atop I\in U}(-1)^{cl(I)}m_I\;
    t^{cl(I)+|I|}},
\end{eqnarray*}
where $[l]=\{1,\ldots,l\}$ and $U\subset 2^{[l]}$ is the ``basis''-set.
However, the conjecture does not include a description of the basis-set $U$. 

Using Algebraic Discrete Morse theory (see \cite{jolwe}), we are able
to specify the basis-set $U$ and prove the conjecture in several cases.
In fact, we give a general conjecture
about the multigraded minimal $A$-free resolution of $k$ over $A$. This 
conjecture implies in these cases an explicit description of the 
multigraded Hilbert and Poincar\'e-Betti series, hence it implies the conjecture by
Charalambous and Reeves.

Section \ref{secmorse} recalls Algebraic Discrete Morse 
theory. For more details and a proof see \cite{jolwe}.

In Section \ref{sec_taylor} we apply Algebraic Discrete Morse theory 
to the Taylor resolution. We define a standard matching 
which we need for the formulation of our conjecture, and we define special 
acyclic matchings for ideals generated in degree two. In particular, we define
matchings (not necessarily acyclic) for Stanley Reisner ideals of order
complexes of a partially ordered set. 

In Section \ref{sec_hilbert} we formulate our conjecture on the 
multigraded minimal resolution of $k$ as an $A$-module and we show that our 
conjecture gives an explicit form of the multigraded Hilbert
and Poincar\'e-Betti series. This generalizes the conjecture by 
Charalambous and Reeves. We say that an algebra $A$ has property \propP~ 
(resp. \propH) if the multigraded Poincar\'e-Betti series (resp. multigraded 
Hilbert series) has the conjectured form.

In Section \ref{sec_koszul} we give a description of the Koszul homology 
$H_\bullet(K^A)$ of the Koszul complex over $A$ with respect to the sequence 
$x_1,\ldots, x_n$ in terms of a standard matching on the Taylor resolution. 
We need this description later in the proof of our conjecture.

In Section \ref{sec_poset} we prove that the Stanley Reisner ring $A=k[\Delta]$, 
where $\Delta=\Delta(P)$ is the order complex of a partially ordered set $P$, 
satisfies property \propP~ and property \propH.

In the first subsection of Section \ref{sec_proof} we prove our conjecture 
for algebras for which $H_\bullet(K^A)$ is an $M$-ring, a notion introduced 
by Fr\"oberg \cite{FR}. Using a theorem of Fr\"oberg, we also prove property
\propP~ for algebras $A=S/\aaf$ for which in addition the minimal free 
resolution of $\aaf$ carries the structure of a differential-graded algebra.
In the second part we prove our conjecture for all Koszul algebras $A$.
Note that this, as a particular case, gives another proof that $A=k[\Delta]$
satisfies property \propP~ and \propH.\\
Finally, we explain why our conjecture makes sense in general.
We generalize the Massey operation in order to get an explicit description
of the Eagon complex. On this complex we define an acyclic matching. If the
resulting Morse complex is minimal, one has to find an isomorphism to the
conjectured complex. We give some ideas on how to construct this isomorphism. 
This construction justifies our conjecture.

Since an algebra is Golod if and only if
\begin{eqnarray*}
  \poinr{k}{A}=\frac{\D\prod_{i=1}^n(1+x_i\;t)}
  {\D 1-t\;\sum_{\beta_{\alpha,i}\neq 0}\beta_{\alpha,i}\ul{x}^\alpha\;t^i},
\end{eqnarray*}
where $\beta_{i,\alpha}:=\dim_k\left(\Tor^S_i(A,k)_\alpha\right)$, we can 
give some applications to the Golod property of monomial rings in the last 
section of this note. We prove, under the assumption of property \propP,~ that 
$A$ is Golod if and only if the first Massey operation is trivial. In 
addition we give, again under the assumption of property \propP, a very 
simple, purely combinatorial condition on the minimal monomial generating 
system $\SM(\aaf)$ which implies Golodness. We conjecture that this is an
equivalence. This would imply that, in the monomial case, Golodness is 
independent of the characteristic of the residue class field $k$.

Recently, Charalambous proved in \cite{chara2} that if 
\[\poinr{k}{A}=\frac{\D\prod_{i=1}^n(1+x_i\;t)}{Q_R(\ul{x},t)}~\mbox{with}~Q_R(\ul{x},t)=\sum\Big(\sum_\alpha c_\alpha \ul{x}^\alpha\Big)t^i,\]
then $\ul{x}^\alpha$ equals to a least common multiple of a subset of 
the minimal monomial generating system $\SM(\aaf)$.
However an explicit form of $Q_R(\ul{x},t)$ in terms of subsets of $\SM(\aaf)$
is still not known.\\
In addition, Charalambous proves a new criterion for generic ideals to be Golod.
In Section \ref{sec_golod} we reprove this criterion using our approach.

In another recent paper, Berglund gives an explicit form of the denominator $Q_R(\ul{x},t)$
in terms of the homology of certain simplicial complexes. Since there seems to be no obvious
connection of the approach taken in \cite{burg} and our approach,
it is an interesting problem to link these two methods.

%
%
\section{Algebraic Discrete Morse Theory}\label{secmorse} 
In this section we recall Algebraic Discrete Morse theory from \cite{jolwe}.

Let $R$ be a ring and $\CRes_\bullet = (C_i, \partial_i)_{i \geq 0}$ be a
chain complex of free $R$-modules $C_i$. We choose a basis
$X=\bigcup_{i=0}^n X_i$ such that $C_i\simeq \Dirsum_{c\in X_i} R\;c$.
>From now on we write the differentials $\partial_i$ with respect to the
basis $X$ in the following form:

$$
\partial_i: \left\{ \begin{array}{lll} C_i & \to & C_{i-1}\\
                                       c & \mapsto & \partial_i(c)=\D{\sum_{c'\in X_{i-1}}}[c:c'] \cdot c'.\\
\end{array} \right.
$$

Given the complex $\CRes_\bullet$ and the basis $X$, we construct a directed,
weighted graph $G(\CRes_\bullet)=(V,E)$. The set of vertices $V$
of $G(\CRes_\bullet)$ is the basis $V=X$ and the set $E$ of (weighted) edges is given
by the rule
\begin{eqnarray*}
(c,c',[c:c'])\in E&:\Leftrightarrow& c\in X_i, c'\in X_{i-1},\mbox{ and }[c:c']\neq 0.
\end{eqnarray*}
We often omit the weight and write $c \rightarrow c'$ to denote an edge in $E$. Also by abuse
of notation we write $e \in G(\CRes_\bullet)$ to indicate that $e$ is an edge in $E$.

\begin{Def} \label{morsedefinition}
A subset $\Mmc\subset E$ of the set of edges is called an
acyclic matching if it satisfies the following three conditions:
\begin{enumerate}
\item (Matching) Each vertex $v\in V$ lies in at most one edge $e\in \Mmc$.
\item (Invertibility) For all edges $(c,c',[c:c'])\in\Mmc$ the weight $[c:c']$ lies
  in the center of $R$ and is a unit in $R$.
\item (Acyclicity) The graph $G_\Mmc(V,E_\Mmc)$ has no directed cycles,
  where $E_\Mmc$ is given by
  \[E_\Mmc:=(E\setminus \Mmc)\cup \left\{\left(c',c,\frac{-1}{[c:c']}\right)
    \mbox{ with }(c,c',[c:c'])\in\Mmc\right\}.\]
\end{enumerate}
\end{Def}

For an acyclic matching $\Mmc$ on the graph $G(\CRes_\bullet) = (V,E)$
we introduce the following notation.

\begin{enumerate}
\item We call a vertex $c \in V$ critical with respect to $\Mmc$
  if $c$ does not lie in an edge $e\in \Mmc$; we write
  $$X^\Mmc_i:=\{c\in X_i~|~c \mbox{~critical~}\}$$ for the set of all critical
  vertices of homological degree $i$.
\item We write $c'\le c$ if $c\in X_i$, $c'\in X_{i-1}$, and $[c:c']\neq 0$.
\item $\on{Path}(c,c')$ is the set of paths from $c$ to $c'$ in the graph $G_\Mmc(\CRes_\bullet)$.
\item The weight $w(p)$ of a path $p=c_1\to\cdots\to c_r\in\on{Path}(c_1,c_r)$ is given by
  \begin{eqnarray*}
    w(c_1\to\ldots\to c_r)&:=&\prod_{i=1}^{r-1} w(c_i\to c_{i+1}),\\
    w(c\to c')&:=&\left\{\begin{array}{rll}
        -\D{\frac{1}{[c:c']}}&,&c\le c',\\ & & \\
        {\D{[c:c']}}&,&c'\le c.\end{array}\right.
  \end{eqnarray*}
\item We write $\Gamma(c,c')=\D{\sum_{p\in \on{Path}(c,c')}}w(p)$ for the sum of
  weights of all paths from $c$ to $c'$.
\end{enumerate}

\noindent Now we are in position to define a new complex $\CRes^\Mmc_\bullet$, which we call
the Morse complex of $\CRes_\bullet$ with respect to $\Mmc$.
The complex $\CRes^\Mmc_\bullet = (C_i^\Mmc,\partial_i^\Mmc)_{i \geq 0}$ is defined by
\[C_i^\Mmc:=\Dirsum_{c\in X^\Mmc_{i}} R\;c,\]

$$\partial_i^\Mmc: \left\{ \begin{array}{ccc} C_i^\Mmc & \to & C_{i-1}^\Mmc\\
                                   c & \mapsto & \D{\sum_{c'\in X^\Mmc_{i-1}}}\Gamma(c,c')c',
\end{array}\right. .$$

\begin{Thm}\label{morse} $\CRes^\Mmc_\bullet$ is a complex of free $R$-modules
and is homotopy-equivalent to the complex $\CRes_\bullet$; in particular, for all $i \geq 0$
\[H_i(\CRes_\bullet) \cong H_i(\CRes^\Mmc_\bullet).\]
The maps defined below give a chain homotopy between $\CRes_\bullet$ and $\CRes^\Mmc_\bullet$:

$$f : \left\{ \begin{array}{lll} \CRes_\bullet & \to     & \CRes^\Mmc_\bullet\\
              c\in X_i             & \mapsto & f(c):=\D{\sum_{c'\in X_i^\Mmc}}\Gamma(c,c')c',
\end{array} \right. $$

$$g : \left\{ \begin{array}{lll} \CRes^\Mmc_\bullet       & \to     & \CRes_\bullet \\
                                c\in X_i^\Mmc & \mapsto & g_i(c):=\D{\sum_{c'\in X_i}}\Gamma(c,c')c'.
\end{array} \right. $$

\end{Thm}

\noindent Sometimes we consider the same construction for matchings which are not 
acyclic. Clearly, Theorem \ref{morse} does not hold anymore for $\CRes_\bullet^{\Mmc}$ if 
$\Mmc$ is not acyclic. In general, there is not even a good definition of the differentials
$\partial^\Mmc$. But for calculating invariants it is sometimes useful to consider 
$\CRes_\bullet^{\Mmc}$ for matchings that are not acyclic. 
In these cases we consider just the vectorspace $\CRes_\bullet^{\Mmc}$.
%
%
%
\section{Algebraic Discrete Morse Theory on the Taylor Resolution}
\label{sec_taylor}
In this section we consider acyclic matchings on the Taylor resolution. 
First, we introduce a {\em standard matching}, which we use 
in later in order to formulate and prove our conjecture.
Then Section \ref{subsec_taylor_grad2} considers the Taylor resolution 
for monomial ideals
which are generated in degree two. The resolutions of those ideals are 
important for the proof of our conjecture in the case where $A$ is 
Koszul (see Section \ref{sec_proof}). 
Next, we give a matching on the Taylor resolution of Stanley 
Reisner ideals of the order complex of a partially ordered set, which we  
use in Section \ref{sec_poset} in order to prove property \propP~ and 
property \propH~ for this type of ideal.\\
Finally, we introduce the (strong) $\gcd$-condition for monomial 
ideals and give a special acyclic matching on the Taylor resolution for this 
type of ideals, which are in connection with the Golod property of monomial 
rings (see Section \ref{sec_golod}). 

\subsection{Standard Matching on the Taylor Resolution}
Let $S=k[x_1,\ldots, x_n]$ be the commutative polynomial ring over a field $k$
of arbitrary characteristic and $\aaf\subg S$ a monomial ideal.

The basis of the Taylor resolution is given by the subsets $I\subset \SM(\aaf)$ 
of the minimal monomial generating system $\SM(\aaf)$ of the ideal $\aaf$. 
For a subset $I\subset \SM(\aaf)$ we denote by $m_I$ the least common 
multiple of the monomials in $I$, $m_I:=\lcm\big(m\in I\big)$.

On this basis we introduce an equivalence relation: We say that two 
monomials $m,n\in I$ with $I\subset \SM(\aaf)$ are equivalent if 
$\gcd(m,n)\neq 1$ and write $m\sim n$. The transitive closure of $\sim$ 
gives us an equivalence relation on each subset $I$. We denote by 
$cl(I):=\# I/\sim$ the number of equivalence classes of $I$.

Based on the Taylor resolution, we define a product by 
$$I\cdot J=\left\{\begin{array}{lcr}0&,&\gcd(m_I,m_J)\neq 1\\I\cup J&,&
\gcd(m_I,m_J)=1.\end{array}\right.$$
Then the number $cl(I)$ counts the factors of $I$ with respect to 
the product defined above.

The aim of this section is to introduce an acyclic matching on the Taylor 
resolution which preserves this product.

We call two subset $I,J\subset \SM(\aaf)$ a {\it matchable pair} and write $I\to J$ 
if $|J|+1=|I|$, $m_J=m_I$, and the differential of the Taylor complex maps 
$I$ to $J$ with coefficient $[I,J]\neq 0$.

Let $I\to J$ be a matchable pair in the Taylor resolution with 
$cl(I)=cl(J)=1$ such that no subset of $J$ is matchable. Then define 
\[\Mmc_{11}:=\{I\dotcup K\to J\dotcup K\mbox{ for each $K$ with 
$\gcd(m_K,m_I)=\gcd(m_K,m_J)=1$}\}.\]
For simplification we write $I\in\Mmc_{11}$ if there exists a subset 
$J$ with $I\to J\in\Mmc_{11}$ or $J\to I\in\Mmc_{11}$. It is clear that 
this is an acyclic matching. Furthermore, the differential changes in each 
homological degree in the same way and for two subsets $I,K$ with 
$\gcd(m_I,m_K)=1$ we have 
$I\dotcup K\in\Mmc_{11}\Longleftrightarrow 
       I\in\Mmc_{11}\mbox{ or }K\in\Mmc_{11}$. 
Because of these facts, we can repeat this matching $\Mmc_{11}$ on the 
resulting Morse complex. This gives us a sequence of acyclic matchings, 
which we denote by $\Mmc_1:=\bigcup_{i\ge 1} \Mmc_{1i}$. If no repetition 
is possible, we reach a resolution with basis 
given by some subsets $I\subset \SM(\aaf)$ with the following property:
If we have a matchable pair $I\to J$ where $I$ has a higher homological 
degree than $J$, then $cl(I)\ge 1$ and $cl(J)\ge 2$. We now construct the 
second sequence:\\
Let $I\to J$ be a matchable pair in the resulting Morse complex 
with $cl(I)=1, cl(J)=2$ such that no subset of $J$ is matchable. Then define 
\[\Mmc_2:=\{I\dotcup K\to J\dotcup K\mbox{ for each $K$ with 
    $\gcd(m_K,m_I)=\gcd(m_K,m_J)=1$}\}.\]
With the same arguments as before this defines an acyclic matching, and a 
repetition is possible. 
The third sequence starts if no repetition of $\Mmc_2$ is possible and is 
given by a matchable pair $I\to J$ in the resulting Morse complex with 
$cl(I)=1, cl(J)=3$ such that no subset of $J$ is matchable. 
Then define 
\[\Mmc_3:=\{I\dotcup K\to J\dotcup K\mbox{ for each $K$ with 
    $\gcd(m_K,m_I)=\gcd(m_K,m_J)=1$}\}.\]
Since every matchable pair is of the form $I\dotcup K\to J\dotcup K$
with $m_I=m_J$, $\gcd(m_I,m_K)=1$, and $cl(I)=1, cl(J)\ge 1$, we finally
reach with this procedure a minimal resolution of the ideal $\aaf$ 
as $S$-module. Let $\Mmc$ be the union of all matchings. 
As before we write $I\in\Mmc$ if there exists a subset $J$ with 
$I\to J\in\Mmc$ or $J\to I\in\Mmc$. Then the minimal resolution 
has a basis given by $\SM(\aaf)\setminus\Mmc$.\\
We give a matching of this type a special name:
\begin{Def}[standard matching] A sequence of matchings $\Mmc:=\bigcup_{i\ge 1}\Mmc_i$ is called a standard 
matching on the Taylor resolution if all the following holds:
\begin{enumerate}
\item $\Mmc$ is graded, i.e. for all edges $I\to J$ in $\Mmc$ we have $m_I = m_J$,
\item $T^\Mmc_\bullet$ is minimal, i.e. for all edges $I\to J$ in $T^\Mmc_\bullet$ we have $m_I\neq m_J$,
\item $\Mmc_i$ is a sequence of acyclic matchings on the Morse complex $T_\bullet^{\Mmc_{<i}}$ 
  ($\Mmc_{<i}:=\bigcup_{j=1}^{i-1}\Mmc_j$, $T_\bullet^{\Mmc_{<1}}=T_\bullet$),
\item for all $I\to J\in\Mmc_i$ we have 
  \begin{eqnarray*}
    cl(J)-cl(I)&=&i-1,\\
    |J|+1&=&|I|,
  \end{eqnarray*}
\item there exists a set $\Bc_i\subset \Mmc_i$ such that 
  \begin{enumerate}
    \item $\D \Mmc_i=\Bc_i\cup~\left\{I\cup K\to J\cup K~\left|~
          \begin{array}{c}\mbox{$K$ with $\gcd(m_I,m_K)=1$}\\
            \mbox{and $I\to J\in \Bc_i$}\end{array}\right.\right\}$ and
    \item for all $I\to J\in \Bc_i$ we have $cl(I)=1$ and $cl(J)=i$.
  \end{enumerate}
\end{enumerate}
\end{Def}

The construction above shows that a standard matching always exists. For a 
standard matching we have two easy properties, which we will need in 
Section \ref{sec_koszul}:
\begin{Lem}\label{canonindep}
Let $\Mmc$ and $\Mmc'$ be two different standard matchings. Then
\begin{enumerate}
\item for all $i\ge 1$ we have $$\D 1+\sum_{I\not\in\Mmc_{<i}}(-1)^{cl(I)}m_I t^{cl(I)+|I|}
  =1+\sum_{I\not\in\Mmc'_{< i}}(-1)^{cl(I)}m_I t^{cl(I)+|I|},$$
\item if $I,J\not\in\Mmc$, $\gcd(m_I,m_J)=1$, and $I\cup J\in \Mmc$, then there exists a set $K$ 
  with $|K|=|I|+|J|+1$, $cl(K)=1$, and $(I\cup J\to K)\in\Mmc$.
\end{enumerate}
\end{Lem}
\begin{proof}
The result follows directly from the definition of a standard matching.
\end{proof}
If the ideal is generated in degree two, every standard matching ends after 
the second sequence: Assume that we have a matchable pair $I\to J$ such that
$cl(I)=1$ and $cl(J)\ge 3$. Then $J$ has at least three subsets $J=J_1\cup J_2\cup J_3$
such that $\gcd(m_{J_i},m_{J_{i'}})=1$, $i,i'=1,2,3$. Since $I$ and $J$ have the same
multidegree and $cl(I)=1$, there would exist a generator $u\in \SM(\aaf)$ such that 
$\gcd(m_{J_i},u)\neq 1$ for $i=1,2,3$. But $u$ is a monomial of degree two, which
makes such a situation impossible.\\
In this case we have
\begin{Lem}\label{secondseq} If every standard matching ends after the second sequence, 
i.e. $\Mmc=\Mmc_1\cup \Mmc_2,$ then
\[\sum_{I\not\in\Mmc_1}(-1)^{cl(I)}m_I t^{cl(I)+|I|}=\sum_{I\not\in\Mmc}(-1)^{cl(I)}m_I t^{cl(I)+|I|}.\]
\end{Lem} 
\begin{proof}
By definition an edge $I\to J$ matched by the second sequence has the 
property $|I|=|J|+1$ and $cl(I)=cl(J)-1$ and $m_I=m_J$. 
Therefore, $$(-1)^{cl(I)}m_I t^{cl(I)+|I|}
= -\left((-1)^{cl(J)}m_J t^{cl(J)+|J|}\right),$$ which proves the assertion.
\end{proof}

\subsection{Resolutions of Monomial Ideals Generated in Degree Two}
\label{subsec_taylor_grad2}
Let $\aaf\subg S$ be a monomial ideal with minimal monomial generating system 
$\SM(\aaf)$ such that for all monomials $m\in \SM(\aaf)$ we have $deg(m)=2$.
We assume, in addition, that $\aaf$ is squarefree. This is no restriction since
via polarization we get similar results for the general case.

First we fix a monomial order $\prec$.
We introduce the following notation:
To each subset $I\subset \SM(\aaf)$ we associate an undirected graph $G_I=(V,E)$ 
on the ground set $V=[n]$, by setting $\{i,j\}\in E$ if the monomial $x_ix_j$ 
lies in $I$. We call a subset $I$ an \nbc-set if the associated graph
$G_I=(V,E)$ contains no broken circuit, i.e. there exists no edge $\{i,j\}$ 
such that
\begin{enumerate}
\item $E\cup\{\{i,j\}\}$ contains a circuit $c$ and
\item $x_ix_j=\max_{\prec}\big\{x_{i'}x_{j'}~\big|~\{i',j'\}\in c\big\}$.
\end{enumerate}
\begin{Prop}\label{nbc-sets}
There exists an acyclic matching $\Mmc_1$ on the Taylor resolution such that 
\begin{enumerate}
\item $\Mmc_1$ is the first sequence of a standard matching,
\item the resulting Morse complex $T_\bullet^{\Mmc_1}$ is a subcomplex of the Taylor 
  resolution and
\item $T_\bullet^{\Mmc_1}$ has a basis indexed by the \nbc-sets. 
\end{enumerate}
\end{Prop}
\begin{proof}
Let $Z$ be a circuit in $T_\bullet$ of maximal cardinality. Let 
$x_ix_j:=\max_{\prec}\{Z\}$. We then define 
\[\Mmc_{1,0}:=\Big\{(Z\cup I)\to((Z\setminus\{x_ix_j\})\cup I)~\Big|~
I\in T_\bullet~\mbox{with}~Z\cap I=\emptyset\Big\}.\]
It is clear that $I$ is an acyclic matching and the resulting Morse complex 
$T^{\Mmc_{1,0}}$ is a subcomplex of the Taylor resolution.\\
Now let $Z_1$ be a maximal circuit in $T^{\Mmc_{1,0}}$ and let 
$x_\nu x_l:=\max_{\prec}\{Z_1\}$. We then define 
\[\Mmc_{1,1}:=\Big\{(Z_1\cup I)\to((Z_1\setminus\{x_\nu x_l\})\cup I)~\Big
|~I\in T^{\Mmc_{1,0}}~\mbox{with}~Z_1\cap I=\emptyset\Big\}.\]
We only have to guarantee that $(Z_1\cup I)\not\in\Mmc_{1,0}$.\\
Assume $(Z_1\cup I)\in\Mmc_{1,0}$. Since 
$(Z_1\setminus\{x_\nu x_l\})\cup I\not\in \Mmc_{1,0}$, we see that $x_\nu x_l\neq x_ix_j$ and 
$x_\nu x_l\in Z$. But then $W:=Z\cup (Z_1\setminus\{x_\nu x_l\})$ is a circuit, which is a 
contradiction to the maximality of $Z$. Therefore, $\Mmc_{1,1}$ is a well defined 
acyclic matching and the resulting Morse complex is a subcomplex of the Taylor 
resolution.\\
If we continue this process, we reach a subcomplex $T^{\Mmc_1}$ of the Taylor 
resolution with a basis indexed by all \nbc-sets. It is clear that 
$\Mmc_1:=\bigcup_i \Mmc_{1,i}$ satisfies all conditions of the first sequence of a 
standard matching. Furthermore, if $I$ is an \nbc-set and $m_I=m_{I\setminus\{m\}}$, 
then it follows that $cl(I)=cl(I\setminus\{m\})-1$ (otherwise we would have a 
circuit). This implies that $\Mmc_1$ is exactly the first sequence of a standard 
matching.
\end{proof}
We denote by $T_{\nbc}$ the resulting Morse complex.
\begin{Cor}Let $\aaf\subg S$ be a monomial ideal generated in degree two. We denote
with $\nbc_i$ the number of \nbc-sets of cardinality $i-1$. Then for the Betti number
of $\aaf$ we have the inequality $\beta_i\le \nbc_i.$
\end{Cor}
\subsection{Resolution of Stanley Reisner Ideals of a Partially Ordered Set}
In this subsection we give a (not acyclic) matching on the subcomplex  
$T_{\nbc}$ in the case where $\aaf=J_{\Delta(P)}$ is the Stanley Reisner 
ideal of the order complex of a partially ordered set $(P,\prec)$. In this 
case $\aaf$ is generated in degree two by monomials $x_ix_j$ where $\{i,j\}$ is 
an antichain in $P$. For simplification we assume that 
$P=[p]=\{1,\ldots, p\}$ and the order $\prec$ preserves the natural order, i.e.
$i\prec j\Rightarrow i<j$, where $<$ is the natural order on the natural 
numbers $\N$. Then the minimal monomial generating system $\SM(\aaf)$ of the 
Stanley Reisner ideal is given by 
\[\SM(\aaf):=\Big\{x_ix_j~\Big|~\mbox{$i<j$ and $i\not\prec j$}\Big\}.\]
Since $\SM(\aaf)$ consists of monomials of degree two, we can work on the 
subcomplex $T_{\nbc}$ of the Taylor resolution, where $T_{\nbc}$ is constructed
with respect to the lexicographic order such that $x_1\succ x_2\succ \ldots \succ x_n$.

First we introduce some notation:
\begin{Def}\label{def_sting}A subset $I\subset \SM(\aaf)$ is called a sting-chain if there 
  exists a sequence of monomials $x_{i_1}x_{i_2},x_{i_2}x_{i_3},\ldots, x_{i_{\nu -1}}x_{i_\nu }\in I$ with 
\begin{enumerate}
    \item $1\le i_1<\ldots<i_\nu \le n$,
    \item $i_1=\min\{j\mbox{ with }x_j\mbox{ divides }\lcm(m_I)\}$,
    \item $i_\nu =\max\{j\mbox{ with }x_j\mbox{ divides }\lcm(m_I)\}$,
    \item for all monomials $x_rx_s\in I$ with $r<s$ exists an index $1\le j\le \nu -1$ such that either 
      \begin{enumerate}
        \item $x_rx_s = x_{i_j}x_{i_{j+1}}$ or
        \item $r=i_{j}$, $s<i_{j+1}$, and $x_sx_{i_{j+1}}\not \in I$ or \label{spikecond1}
        \item $r>i_{j}$, $s=i_{j+1}$, and $i_j\prec r$ (i.e. $x_{i_j}x_r\not \in \aaf$).\label{spikecond2}
      \end{enumerate}    
  \end{enumerate}
\end{Def}
Let $\Bc$ be the set of all chains of sting-chains:
\[\Bc:=\left\{(I_1,\ldots, I_l)\left|
\begin{array}{c}I_j\mbox{ sting-chain for all $j=1,\ldots, l$ and }\\
      \max(I_j)<\min(I_{j+1})\mbox{ for all $j=1,\ldots, l-1$}\end{array}\right.\right\},\]
where 
\[\begin{array}{l}\max(I):=\max\{i\mid x_i~\mbox{divides}~\lcm(m_I)\}\\   
\min(I):=\min\{i\mid x_i~\mbox{divides}~\lcm(m_I)\}.\end{array}\]

Note that a sting-chain is not necessarily an {\nbc}-set. For example, the set 
$\{x_ix_l, x_\nu x_l, x_jx_l\}$ with $i<\nu <j<l$ is a sting-chain, if 
$x_ix_\nu , x_ix_j\not\in\aaf$, but it contains a broken circuit if $x_\nu x_j\in\aaf$.
But with an identification of those sets we get the following Proposition:
\begin{Prop}\label{spikechain}There exists a matching $\Mmc_2$ (not necessary acyclic) on the complex 
$T_{\nbc}$ such that 
\begin{enumerate}
\item there exists a bijection between the sets $I\in T_{\nbc}^{\Mmc_2}$ and the 
  chains of sting-chains $I\in\Bc$,
\item for $I\to I'\in\Mmc_2$ we have 
  \begin{enumerate}
    \item $\lcm(m_I)=\lcm(m_{I'})$ and
    \item $cl(I)=cl(I')-1$ and $|I|=|I'|+1$.
  \end{enumerate}
\end{enumerate} 
\end{Prop}
\begin{proof}
For a set $I\in T_{\nbc}\setminus\mathcal{B}$ let $x_ix_\nu x_jx_l$ be the maximal monomial
with respect to the lexicographic order such that $i<\nu <j<l$ and at least one of the 
following conditions is satisfied:
\begin{enumerate}
\item $x_ix_j, x_\nu x_l\in I$ and $x_ix_l\not\in I$,
\item $x_ix_l, x_\nu x_j\in I$.
\end{enumerate}
{\bf Case $x_ix_j, x_\nu x_l\in I$:} Because of the transitivity of the order 
  $\prec$ on $P$ we have either $x_ix_\nu \in\aaf$ or $x_\nu x_j\in\aaf$. 
  \begin{enumerate}
    \item[$\triangleright$] Assume $x_ix_\nu \in\aaf$. Since $x_ix_\nu x_jx_l$ is the 
      maximal monomial satisfying one of the conditions above, it follows that if 
      $I\cup\{x_ix_\nu \}$ contains a broken circuit, then $I\setminus\{x_ix_\nu \}$ 
      contains a broken circuit as well. We set  
      \[\Big((I\setminus\{x_ix_\nu \})\dotcup J\Big)\to \Big((I\cup\{x_ix_\nu \})\dotcup J\Big)\in\Mmc_2
      \mbox{\hspace{2cm}}\]
      \[\mbox{ for all $J$ with $\gcd(\lcm(I),\lcm(J))=1$}.\]
    \item[$\triangleright$] If $x_ix_\nu \not\in\aaf$, then $x_\nu x_j\in\aaf$. Again, we 
      have that if $I\cup\{x_\nu x_j\}$ contains a broken circuit, then 
      $I\setminus\{x_\nu x_j\}$ contains a broken circuit as well. In this case we set
      \[\Big((I\setminus\{x_\nu x_j\})\dotcup J\Big)\to \Big((I\cup\{x_\nu x_j\})\dotcup J\Big)\in\Mmc_2
      \mbox{\hspace{2cm}}\]
      \[\mbox{~\hspace{2cm}for all $J$ with $\gcd(\lcm(I),\lcm(J))=1$}.\]
  \end{enumerate}
{\bf Case $x_ix_l, x_\nu x_j\in I$:} Again, the transitivity implies $x_ix_\nu \in\aaf$ or 
  $x_\nu x_l\in\aaf$ and $x_ix_j\in\aaf$ or $x_jx_l\in\aaf$: 
  \begin{enumerate}
    \item[$\triangleright$]Assume $x_ix_\nu \in\aaf$. As above we have that if 
      $I\cup\{x_ix_\nu \}$ contains a broken circuit, then $I\setminus\{x_ix_\nu \}$ 
      contains a broken circuit as well. We set 
      \[\Big((I\setminus\{x_ix_\nu \})\dotcup J\Big)\to \Big((I\cup\{x_ix_\nu \})\dotcup J\Big)\in\Mmc_2
      \mbox{\hspace{2cm}}\]
      \[\mbox{~\hspace{2cm}for all $J$ with $\gcd(\lcm(I),\lcm(J))=1$}.\]
    \item[$\triangleright$]If $x_ix_\nu \not\in\aaf$, then $x_\nu x_l\in\aaf$. Assume 
      $x_ix_j\in\aaf$. Then again we have that if $I\cup\{x_ix_j\}$ contains a 
      broken circuit, then $I\setminus\{x_ix_j\}$ also contains a broken circuit. 
      In this case we set
      \[\Big((I\cup\{x_ix_j\})\dotcup J\Big)\to\Big((I\setminus\{x_ix_j\})\dotcup J\Big)\in\Mmc_2
      \mbox{\hspace{2cm}}\]
      \[\mbox{~\hspace{2cm}for all $J$ with $\gcd(\lcm(I),\lcm(J))=1$}.\]
    \item[$\triangleright$]Now assume $x_ix_\nu , x_ix_j\not\in\aaf$, then 
      $x_\nu x_l, x_jx_l\in\aaf$. Assume further that $x_jx_l\not\in I$. Then we set
      \[\Big((I\cup\{x_\nu x_l\})\dotcup J\Big)\to\Big((I\setminus\{x_\nu x_l\})\dotcup J\Big)\in\Mmc_2
      \mbox{\hspace{2cm}}\]
      \[\mbox{~\hspace{2cm}for all $J$ with $\gcd(\lcm(I),\lcm(J))=1$}.\]
    \item[$\triangleright$]Finally, we have to discuss the case 
      $x_ix_\nu , x_ix_j\not\in\aaf$ and $x_jx_l\in I$. Then the set $I$ cannot be 
      matched because adding $x_\nu x_l$ would give a circuit and by removing $x_jx_l$ 
      we get a set which is already matched. We identify these sets with the 
      sets containing $x_ix_l, x_\nu x_l, x_jx_l$ instead of $x_ix_l, x_\nu x_j, x_jx_l$. 
      Therefore, this case gives us all sets which are sting-chains 
      but not \nbc-sets.
  \end{enumerate}

With the identification we can say that an \nbc-set $I\not\in\Mmc$ satisfies the 
following two properties, which are exactly the properties of $I\in\Bc$:
\begin{enumerate}
\item If there exist $i<\nu <j<l$ such that $x_ix_j,x_\nu x_l\in I$, then $x_ix_l\in I$ and 
  $x_\nu x_j,x_jx_l\not\in I$ and $x_ix_\nu \not\in\aaf$.
\item There exist no $i<\nu <j<l$ such that $x_ix_l,x_\nu x_j\in I$.
\end{enumerate}
\end{proof}
Note that $T^{\Mmc_2}$ is not a resolution (not even a complex), but we need it because of the 
following corollary, which will be important in Section \ref{sec_poset}.
\begin{Cor}\label{sum_koszul}Let $\aaf$ be a monomial ideal generated in 
degree two and $\Mmc=\Mmc_1\cup \Mmc_2$ a standard matching on the Taylor 
resolution. With the notation above we get:
\begin{equation}
  \begin{array}{lll}
    \D \sum_{I\not\in\Mmc_1}(-1)^{cl(I)}m_I t^{cl(I)+|I|}
    &=&\D \sum_{I\not\in\Mmc}(-1)^{cl(I)}m_I t^{cl(I)+|I|}\\
    &=&\D \sum_{I~\sst \mbox{\nbc-set}}(-1)^{cl(I)}m_I t^{cl(I)+|I|}.
  \end{array}\label{nbc-equa}
\end{equation}
If $\aaf$ is the Stanley Reisner ideal of the order complex of a partially ordered set $P$, then
\begin{equation}
(\ref{nbc-equa})=\sum_{I\not\in\Bc}(-1)^{cl(I)}m_I t^{cl(I)+|I|}.\label{spike-eq}
\end{equation}
\end{Cor}
\begin{proof}
Lemma \ref{secondseq} implies the first equality and the second equality 
follows by Proposition \ref{nbc-sets}. If $\aaf$ is the Stanley Reisner ideal of the 
order complex of a partially ordered set $P$, then Proposition \ref{spikechain}
together with the proof of Lemma \ref{secondseq} imply Equation
(\ref{spike-eq}). 
\end{proof}

\subsection{The $\gcd$-Condition}
In this subsection we introduce the $\gcd$-condition. Let $\aaf\subg S$ be a 
monomial ideal in the commutative polynomial ring and $\SM(\aaf)$ a minimal 
monomial generating system.
\begin{Def}[$\gcd$-condition]
\begin{enumerate}
\item We say that $\aaf$ satisfies the $\gcd$-con\-dition, if for any two monomials 
  $m, n\in \SM(\aaf)$ with $\gcd(m,n)=1$ there exists a monomial $m,n\neq u\in \SM(\aaf)$ 
  with $u\mid\lcm(m,n)$;
\item We say that $\aaf$ satisfies the strong $\gcd$-condition if there exists a 
  linear order $\prec$ on $\SM(\aaf)$ such that for any two monomials 
  $m\prec n\in \SM(\aaf)$ with $\gcd(m,n)=1$ there exists a monomial $m,n\neq u\in \SM(\aaf)$ 
  with $m\prec u$ and $u\mid\lcm(m,n)$.
\end{enumerate}
\end{Def}

\begin{Ex}Let $\aaf=\langle x_1x_2,x_2x_3,x_3x_4,x_4x_5,x_1x_5\rangle$ be the 
  Stanley Reisner ideal of the triangulation of the $5$-gon. Then $\aaf$ 
  satisfies the $\gcd$-condition, but not the strong $\gcd$-condition.
\end{Ex}

\begin{Prop}\label{gcd-complex} Let $\aaf$ be a monomial ideal which satisfies the 
  strong $\gcd$-condition. Then there exists an acyclic matching $\Mmc$ on the Taylor 
  resolution such that for all $\SM(\aaf)\supset I\not\in\Mmc$ we have $cl(I)=1$. 
  We call the resulting Morse complex $T_{\bf \gcd}$.
\end{Prop}
\begin{proof}
Assume $\SM(\aaf)=\{m_1\prec m_2\prec\ldots\prec m_l\}$. We start with $m_1$. 
Let $m_{i_0}\in \SM(\aaf)$ be the smallest monomial such that $\gcd(m_1,m_{i_0})=1$. 
Then there exists a monomial 
$m_1\prec u_0\in \SM(\aaf)$ with $u_0\mid \lcm(m_1,m_{i_0})$. Then we define
\[\Mmc_0:=\Big\{\Big(\{m_1,m_{i_0},u_0\}\cup I\Big)\to \Big(\{m_1,m_{i_0}\}\cup I\Big)~\Big|~I\subset \SM(\aaf)\Big\}.\]
It is clear that this is an acyclic matching and that the Morse complex 
$T_\bullet^{\Mmc_0}$ is a subcomplex of the Taylor resolution.\\
Now let $m_{i_1}$ be the smallest monomial $\neq m_{i_0}$ such that $\gcd(m_1,m_{i_1})=1$.
Then there exists a monomial $m_1\prec u_1\in \SM(\aaf)$ with $u_1\mid \lcm(m_1,m_{i_1})$ and 
we define
\[\Mmc_1:=\Big\{\Big(\{m_1,m_{i_1},u_1\}\cup I\Big)\to \Big(\{m_1,m_{i_1}\}\cup I\Big)~\Big|~I\subset \SM(\aaf)\Big\}.\]
Again, it is straightforward to prove that $\Mmc_1$ is an acyclic matching on 
$T^{\Mmc_0}$ and that the Morse complex is a subcomplex of the Taylor 
resolution. We repeat this process for all $m_1\prec m_i$ with 
$\gcd(m_1,m_i)=1$ and we reach a subcomplex $T^{\Mmc_{m_1}}$, 
$\Mmc_{m_1}=\bigcup_i\Mmc_i$, of the Taylor resolution which satisfies the 
following condition: For all remaining subsets 
$I\subset \SM(\aaf)\setminus\Mmc_{m_1}$ we have:
\begin{enumerate}
\item $m_1\in I\Rightarrow cl(I)=1$,\\\label{welldefined}
\item $m_1\not\in I\Rightarrow cl(I)\ge 1$.
\end{enumerate}
We repeat now this process with the monomial $m_2$. Here we have to guarantee that for 
a set $\{m_2,m_i\}\cup I$ the corresponding set $\{m_2,m_i,u_i\}\cup I$, with 
$\gcd(m_2,m_i)=1$ and $m_2\prec u_i$ and $u_i\mid\lcm(m_2,m_i)$, is not matched by the 
first sequence $\Mmc_{m_1}$.
Since all sets $J\in\Mmc_{m_1}$ satisfy $m_1\in J$, this would be the case if either 
$u_i=m_1$ or $m_1\in I$. The first case is impossible since $m_1\prec m_2\prec u_i$. 
In second case we have $cl\big(\{m_2,m_i\}\cup I\big)=1$. We define:
\[\Mmc_2:=\left\{\Big(\{m_2,m_i,u_2\}\cup I\Big)\to\Big(\{m_2,m_i\}\cup I\Big)
~\left|~\begin{array}{c}I\subset \SM(\aaf)\setminus\Mmc_{m_1}\\
\mbox{and}~cl\big(\{m_2,m_i\}\cup I\big)\ge 2\end{array}\right.\right\}.\]
Condition (\ref{welldefined}) implies then that $\Mmc_2$ is a well defined sequence of 
acyclic matchings. Since we make this restriction, the resulting Morse complex is not 
anymore a subcomplex of the Taylor resolution, but we have still the following fact:
For all remaining subsets $I\subset \SM(\aaf)\setminus\big(\Mmc_{m_1}\cup\Mmc_{m_2}\big)$ 
we have:
\begin{enumerate}
\item $m_1\in I\Rightarrow cl(I)=1$,\\
\item $m_2\in I\Rightarrow cl(I)=1$,\\
\item $m_1,m_2\not\in I\Rightarrow cl(I)\ge 1$.
\end{enumerate}
We apply this process to all monomials. Then we finally reach a complex with the 
desired properties.
\end{proof}

%
%
%
\section{The Multigraded Hilbert and Poincar\'e-Betti Series}\label{sec_hilbert}
Let $\aaf\subg S$ be a monomial ideal and $\Mmc=\Mmc_1\cup\bigcup_{i\ge 2}
\Mmc_i$ a standard matching on the Taylor resolution.
We introduce a new non-commutative polynomial ring $\tilde{R}$, defined by
\[\tilde{R}:=k\langle Y_I~\mbox{for}~\SM(\aaf)\supset I\not\in\Mmc_1~
\mbox{and}~cl(I)=1\rangle.\]
On this ring, we define three gradings:
\begin{eqnarray*}
|Y_I|&:=&|I|+1,\\
\deg(Y_I)&:=&\alpha,~\mbox{with}~\ul{x}^\alpha=m_I,\\
\deg_t(Y_I)&:=&||\alpha||,~\mbox{with}~\ul{x}^\alpha=m_I,
\end{eqnarray*}
where $||\alpha||=\sum_i \alpha_i$ is the absolute value of $\alpha$.
This makes $\tilde{R}$ into a multigraded ring: 
$$\tilde{R}=\bigoplus_{\alpha\in\N^n}\bigoplus_{i\ge 0}~\tilde{R}_{i,\alpha}$$
with $\tilde{R}_{i,\alpha}:=\big\{u\in \tilde{R}~\big|~\deg(u)=
\alpha~\mbox{and}~|u|=i\big\}.$

Let $[Y_I,Y_J]:=Y_IY_J-(-1)^{|Y_I||Y_J|}Y_JY_I$ be the graded commutator 
of $Y_I$ and $Y_J$. We define the following multigraded two-side ideal 
$$\rr:=\langle [Y_I,Y_J]~\mbox{for}~\gcd(m_I,m_J)=1\rangle,$$
and set \[R:=\tilde{R}/\rr.\]
Let $\D \Hilb_R(\ul{x},t,z):=\sum_{\alpha\in\N^n}
\sum_{i\ge 0}\dim_k(R_{i,\alpha})~\ul{x}^\alpha~t^{||\alpha||}~z^i$
be the multigraded Hilbert series of $R$. We have the following fact:
\begin{Prop}The multigraded Hilbert series $\Hilb_R(\ul{x},t,z)$ of $R$ is 
given by
\[\Hilb_R(\ul{x},t,z)=\frac{1}{\D 1+\sum_{I\subset \SM(\aaf)\atop I\not\in\Mmc_1}
  (-1)^{cl(I)}~m_I~t^{m_I}~z^{cl(I)+|I|}},\]
where $t^{m_I}:=t^\alpha$ with $\ul{x}^\alpha=m_I$.
\end{Prop}
\begin{proof}
In \cite{foata}, Cartier and Foata prove that the Hilbert series of an 
arbitrary non-commutative polynomial ring divided by an ideal, which 
is generated by some (graded) commutators, is given by 
\[\Hilb_R(\ul{x},t,z):=\frac{\D 1}{\D 1+\sum_{F}(-1)^{|F|}~\ul{x}^{\deg(y_F)}~
  t^{\deg_t(y_F)}~z^{|Y_F|}},\]
where $F\subset\{Y_I\mbox{ with }I\not\in\Mmc_1, cl(I)=1\}$ is a commutative 
part (i.e. $Y_IY_J=(-1)^{|J||I|}Y_JY_I$ for all $Y_I,Y_J\in F$) and 
$Y_F=\prod_{Y_I\in F}Y_I$.\\
Therefore, we only have to calculate the commutative parts. Since $\rr$ is 
generated by the relations 
$Y_IY_J=(-1)^{|J||I|}Y_JY_I$, if $\gcd(m_I,m_J)=1$, we see that 
the commutative parts are given by 
$$F:=\Big\{Y_{I_{i_1}},\ldots, Y_{I_{i_r}}~\Big|~
\gcd(m_{I_{i_j}},m_{I_{i_{j'}}})=1~\mbox{for all $j\neq j'$}\Big\}.$$
But the fact that $Y_{I_{i_1}},\ldots, Y_{I_{i_r}}$ is a commutative part is 
equivalent to $I_{i_1}\cup\ldots\cup I_{i_r}\not\in\Mmc_1$. Therefore, we 
can identify the commutative parts $F$ with the elements 
$I\not\in\Mmc_1$ and sum over all $I\not\in\Mmc_1$.
It is clear that the cardinality of a commutative part equals to the number 
$cl(I)$. If $I=I_1\dotcup\ldots\dotcup I_r$ with $cl(I_j)=1$ is a 
commutative part, it follows that $Y_I=Y_{I_1}\cdots Y_{I_r}$, which implies 
the exponents of $t,z,\ul{x}$.
\end{proof}
%
%
%
%
We formulate the following conjecture:
\begin{Conj}\label{poin-conj}
Let $F_\bullet$ be a multigraded minimal $A$-free resolution of $k$ as an 
$A$-module with 
$F_i:=\bigoplus_{\alpha\in\N^n}A(-\alpha)^{\beta_{i,\alpha}}$ for $i\ge 0$.
Then we have the following isomorphism as $k$-vectorspaces:
\[F_i\iso \bigoplus_{J\subset\{1,\ldots,n\}\atop |J|=l}
\bigoplus_{u\in R\atop |u|=i-l}~A\big(-(\alpha_J+\deg(u))\big),\]
where $\alpha_J$ is the characteristic vector of $J$, defined by
\[(\alpha_J)_i=\left\{\begin{array}{lcr}0&,&i\not\in J,\\1&,&i\in J.\end{array}
\right.\]
\end{Conj}
%
%
This conjecture gives a precise formulation of the conjecture by Charalambous and Reeves
on the multigraded Poincar\'e-Betti series.
In addition, we get an explicit form of the multigraded Hilbert series of $S/\aaf$ for
monomial ideals $\aaf$.
\begin{Prop}Let $A=S/\aaf$ be the quotient 
of the commutative polynomial ring by a monomial ideal $\aaf$, and let 
$\Mmc:=\Mmc_1\cup\bigcup_{i\ge 2}\Mmc_i$ be a standard matching on the Taylor 
resolution. If Conjecture \ref{poin-conj} holds, then the multigraded Poincar\'e-Betti and 
Hilbert series have the following form:
\begin{eqnarray}
  \poinr{k}{A}&=&\prod_{i=1}^n(1+x_i~t)~\Hilb_R(\ul{x},1,t)\label{poinseries}\\
  \nonumber &=&\frac{\D\prod_{i=1}^n(1+x_i~t)}
  {\D 1+\sum_{I\subset \SM(\aaf)\atop I\not\in \Mmc_1}(-1)^{cl(I)}~m_I~
    t^{cl(I)+|I|}},
\end{eqnarray}
\begin{eqnarray}
  \hilb{A}&=&\Big(\prod_{i=1}^n(1-x_i~t)~
  \Hilb_R(\ul{x},t,-1)\Big)^{-1}\label{hilbseries}\\
  \nonumber &=&\frac{\D 1+\sum_{I\subset \SM(\aaf)\atop I\not\in \Mmc_1}
    (-1)^{|I|}~m_I~t^{m_I}}{\D\prod_{i=1}^n(1-x_i~t)}.
\end{eqnarray}
\end{Prop}
Note that Equation (\ref{poinseries}) is a reformulation of the conjecture by 
Charalambous and Reeves. 
\begin{proof}
The form of the Poincar\'e-Betti series follows directly from the conjecture, by 
counting basis elements of $F_i$.\\
For the Hilbert series we consider the complex $F_\bullet\to k\to 0$, which 
is exact since $F_\bullet$ is a minimal free resolution of $k$. Since the 
Hilbert series of $k$ is $1$, the Euler characteristic implies:
\[\sum_{i\ge 0}(-1)^i~\hilb{F_i}=1.\]
Conjecture \ref{poin-conj} implies 
\[\hilb{F_i}=\sum_{J\subset\{1,\ldots,n\}\atop |J|=l}\sum_{u\in R\atop |u|=i-l}
    \ul{x}^{\alpha_J}~t^{|J|}~~\ul{x}^{\deg(u)}~t^{\deg_t(u)}~~\hilb{A}.\]
The Cauchy product finally implies:
\begin{eqnarray*}
\sum_{i\ge 0}(-1)^i\hilb{F_i}&=&\hilb{A}~\sum_{i\ge 0}
\sum_{J\subset\{1,\ldots,n\}\atop |J|=l}(-1)^l~\ul{x}^{\alpha_J}~t^{|J|}\\
&&\sum_{u\in R\atop |u|=i-l}(-1)^{i-l}~\ul{x}^{\deg(u)}~t^{\deg_t(u)}\\
&=&\hilb{A}~\left(\sum_{J\subset\{1,\ldots,n\}}\ul{x}^{\alpha_J}~(-t)^{|J|}
\right)\\
&&\left(\sum_{u\in R}\ul{x}^{\deg(u)}~t^{\deg_t(u)}~(-1)^{|u|}\right)\\
&=&\hilb{A}~\prod_{i=1}^n(1-t~x_i)~\Hilb_R(\ul{x},t,-1).
\end{eqnarray*}
\end{proof}
It is known that if $A$ is Koszul, then $\hilb{A}=1/\poinm{k}{A}$. In our case,
this means:
\begin{Prop}
If $A$ is Koszul, then $\Hilb_R(\ul{x},t,-1)=\Hilb_R(\ul{x},1,-t)$.
\end{Prop}
\begin{proof}
In the monomial case, the Koszul property is equivalent to the fact that $\aaf$ 
is generated in degree two. We prove that a subset $I\in \SM(\aaf)$ which 
is not matched by $\Mmc_1$ satisfies $cl(I)+|I|=\deg_t(Y_I)$.
It is clear that this proves the assertion.\\
It is enough to prove it for subsets $I\subset \SM(\aaf)$ with $cl(I)=1$. 
Let $m_I=\ul{x}^\alpha$ be the least common multiple of the generators in $I$. 
Since all generators have degree two, it follows 
$||\alpha||\le 2+|I|-1=|I|+1=|I|+cl(I)$. 
Since $\Tor^S_i(S/\aaf,k)_i=0$, we get $||\alpha||=|I|+1=|I|+cl(I)$.
\end{proof}
We introduce some notation for rings $A$ satisfying the consequences of Conjecture 
\ref{poin-conj}.
\begin{Def} We say that $A$ has property 
\begin{enumerate}
\item[\propP]if $\poinr{k}{A}=\prod_{i=1}^n(1+x_i~t)~\Hilb_R(\ul{x},1,t)$ and 
has property 
\item[\propH]if $\hilb{A}=\Big(\prod_{i=1}^n(1-x_i~t)~\Hilb_R(\ul{x},t,-1)
\Big)^{-1}$.
\end{enumerate}
\end{Def}

%
%
%
\section{The Homology of the Koszul Complex $K^A$}\label{sec_koszul}
Let $\mathcal {M}$ be a standard matching on the Taylor resolution of $\aaf$. 
The basis of the $k$-vectorspace $T_\bullet^\Mmc\otimes_S k$ is then given by 
the sets $I\subset \SM(\aaf)$ with $I\not\in\Mmc$.

We denote with $K^A_\bullet$ the Koszul complex of $A$ with respect to the sequence 
$x_1,\ldots, x_n$, i.e.
\[K_i:=\bigoplus_{\{j_1<\ldots<j_i\}}A~e_{\{j_1<\ldots<j_i\}}\]
with differential 
\[\partial_i:\left\{\begin{array}{ccl}K_i&\to& K_{i-1}\\
    e_{\{j_1<\ldots<j_i\}}&\mapsto&
    \sum_{l=1}^i (-1)^{l+1}~x_{j_l}~e_{\{j_1<\ldots<j_{l-1}<j_{l+1}<\ldots j_i\}}\end{array}\right.\]
We denote further by $Z(K_\bullet)$ (resp. $B(K_\bullet)$) the set of cycles (resp. boundaries)
of the complex $K_\bullet$. Finally, we denote with $H(K_\bullet)$ the homology of the Koszul 
complex.
\begin{Prop}\label{homo}
If $\Mmc$ is a standard matching, then there exists a homogeneous homomorphism 
\[\phi:\left\{\begin{array}{ccc}T^\Mmc_\bullet\otimes_S k&\to&K^A_\bullet\\
I&\mapsto&\phi(I)\end{array}\right.\]
such that for all $I,J\not\in\Mmc$ with $\gcd(m_I,m_J)=1$ we have 
\begin{enumerate}
 \item $\phi(I)$ is a cycle,
 \item $\phi(I)\phi(J) = \phi(I\cup J)$ if $I\cup J\not\in\Mmc$, 
 \item if $I\cup J\in\Mmc$,
   \[\phi(I)\phi(J)=\partial(c)+
       \sum_{\begin{array}{c}\sst L\not\in\Mmc\\\sst cl(L)\ge cl(I)+cl(J)\end{array}}a_L \phi(L)
     \mbox{~\hspace{5mm}for some $a_L\in k$},\]
     for some $c\in K^A_\bullet$.
\end{enumerate}
Note that $\phi(I)\phi(J)\in B(K_\bullet)$ might happen if all coefficients $a_L$ are zero.
\end{Prop}
\begin{proof} We consider the following double complex:
\[\begin{array}{ccccccccccc}
    &    & 0                        &   &      &   &   0                       &   &   0           &   &\\
    &    &\uparrow                  &   &      &   &\uparrow                   &   & \uparrow      &   &\\
  0 &\to &T_n^\Mmc\otimes_S k&\to&\ldots&\to& T_0^\Mmc\otimes_S k&\to& S/I\otimes_S k&\to& 0\\
    &    &\uparrow                  &   &      &   &\uparrow                   &   & \uparrow      &   &\\  
  0 &\to &T_n^\Mmc\otimes_S K^S_0&\to&\ldots&\to& T_0^\Mmc\otimes_S K^S_0&\to& S/I\otimes_S K^S_0&\to& 0\\
  &    &\uparrow                  &   &      &   &\uparrow                   &   & \uparrow      &   &\\
  0 &\to & \vdots                 &\to&\ldots&\to&\vdots                     &\to& \vdots        &\to& 0\\
  &    &\uparrow                  &   &      &   &\uparrow                   &   & \uparrow      &   &\\
  0 &\to &T_n^\Mmc\otimes_S K^S_n&\to&\ldots&\to& T_0^\Mmc\otimes_S K^S_n&\to& S/I\otimes_S K^S_n&\to& 0\\
  &    &\uparrow                  &   &      &   &\uparrow                   &   & \uparrow      &   &\\
  &    & 0                        &   &      &   &   0                       &   &   0         &   &
\end{array}\]
Since every row and every column, except the first row and the right column, are exact, we 
get by diagram chasing a homogeneous homomorphism 
\[\phi:\left\{\begin{array}{ccc}T^\Mmc_\bullet\otimes_S k&\to&K_\bullet\\
I&\mapsto&\phi(I).\end{array}\right.\]
By construction it is clear that $\phi(I)$ is a cycle.
The second condition of a standard matching is: if $(I\to J)\in \Mmc$, then 
$(I\cup K\to J\cup K)\in\Mmc$ for all $K$ with $\gcd(m_K,m_I)=1$.
This condition implies that one can chose the homomorphism $\phi$ such that 
$\phi(I)\phi(J) = \phi(I\cup J)$ if $I\cup J\not\in\Mmc$.\\
Now let $I\cup J\in\Mmc$. Since $I,J\not\in\Mmc$, it follows from the 
standard matching that $I\cup J$ is matched with a set $\hat I$ of higher homological 
degree. We now consider $\Mmc':=\Mmc\setminus\{\hat I\to I\cup J\}$. We then have 
\[0=\partial^{\Mmc'}\partial^{\Mmc'}(\hat I).\]
Hence we get:
\[\partial^{\Mmc'}(I\cup J)=\sum_{L\not\in\Mmc} a_L \partial^\Mmc(L).\]
Since we take the tensor product $\otimes_S k$ with $k$, all summands with $a_L\not\in k$ 
cancel out. Hence $\phi(I)\phi(J)\in B(K^A_\bullet)$ or, again with diagram chasing:
\[\phi(I)\phi(J)-\sum_{\sst L\not\in\Mmc\atop\sst cl(L)\ge cl(I)+cl(J)}a_L \phi(L)\in B(K_\bullet^A).\]
>From the construction of the standard matching it follows, in addition, 
that $cl(L)\ge cl(I)+cl(J)$ (otherwise $L$ would have been matched before).
\end{proof}

\medskip

We define the following new $k$-algebra:\\
For each $I\not\in\Mmc$ with $cl(I)=1$ we define one indeterminate $Y_I$ with 
total degree $\deg_t(Y_I):=|I|$ and multidegree $\deg_m(Y_I):=x^\alpha$, 
if $x^\alpha=m_I$. Let $R':=k(Y_I, I\not\in\Mmc, cl(I)=1)/\rr'$ be the 
quotient algebra of the graded commutative polynomial ring $k(Y_I, I\not\in\Mmc, cl(I)=1)$ 
(i.e. $Y_IY_J=(-1)^{|I||J|}Y_JY_I$) and the multigraded ideal $\rr'$ 
that is generated by the relations given by Proposition \ref{homo}, i.e.:
\begin{enumerate}
\item $Y_IY_J=0$ if $\gcd(m_I,m_J)\neq 1$,
\item $Y_{I_{i_1}}\cdots Y_{I_{i_r}}=\sum a_L Y_L$ if 
  $\phi(I_{i_1})\cdots \phi(I_{i_r})=\sum a_L \phi(L)+boundary$,
\item $Y_{I_{i_1}}\cdots Y_{I_{i_r}}=0$ if $[\phi(I_{i_1})\cdots \phi(I_{i_r})]=0$.
\end{enumerate} 
\begin{Thm}\label{mring}
If $\Mmc$ is a standard matching, then $R'$ is isomorphic to $H(K_\bullet)$.
\end{Thm}
\begin{proof} The isomorphism is given by Proposition \ref{homo}. 
We only have to prove that $[\phi(I)][\phi(J)]=0$ if $\gcd(m_I,m_J)\neq 1$. 
This follows from the next lemma and the next corollary.
\end{proof}
\begin{Lem}\label{einfach}
Let $c=\sum_I \alpha_I \frac{m}{x_I}\,e_I$ be a homogeneous cycle with multidegree $\deg(c)=m$. 
We fix an $x_0\mid m$. Then there exists a 
cycle $c'=\sum_{I'}\alpha_{I'}\frac{m}{x_{I'}}e_{I'}$, homologous to $c$, 
such that $x_0\mid x_{I'}$ for all $I'$.
\end{Lem}
\begin{proof} Let $I$ be an index set such that $\alpha_I\neq 0$ in the expansion of $c$ with 
$x_0\not| ~x_I$. Then
\begin{eqnarray}
\frac{m}{x_I}e_I = \sum_{i\in I}(-1)^{pos(i)+1}\frac{m\,x_i}{x_0\,x_I}
   e_{x_0}\wedge e_{I\setminus\{i\}}+\partial\left(\frac{m_I}{x_0\,x_I}e_{x_0}\wedge e_I\right).\label{ers}
\end{eqnarray}
If we replace each index set $I$ with respect to (\ref{ers}), we finally reach a cycle $c'$ 
with the desired property. By construction there exists an element $d$ with 
$c-c'=\partial(d)\in B(K_\bullet)$. 
\end{proof}
\begin{Cor}\label{immer}
Let $c_1,c_2$ be two homogeneous cycles with multidegrees $\deg(c_1)=m$ and $\deg(c_2)=n$. 
If $\gcd(m,n)\neq 1$, we have $[c_1][c_2]=0$.
\end{Cor}
\begin{proof} Let $c_1:=\sum_I \alpha_I\frac{m}{x_I}e_I$ and  
$c_2:=\sum_J \beta_J\frac{n}{x_J}e_J$ with $\gcd(m,n)\neq 1$. 
We fix a $j\in\supp(\gcd(m,n))$. By Lemma \ref{einfach} we can assume 
that $j\in I\cap J$ for all $I,J$. This implies $[c_1][c_2]=0$.
\end{proof}
\begin{Cor}$H(K_\bullet)$ is generated by $I\not\in\Mmc$ with $cl(I)=1$.
\end{Cor}

%
%
%
\section{Hilbert and Poincar\'e-Betti Series of the Algebra $A=k[\Delta]$}
\label{sec_poset}
In this section we prove property \propP~ and \propH~ for $A=S/\aaf$ where $\aaf=I_{\Delta(P)}$ 
is the Stanley Reisner ideal of the order complex $\Delta(P)$ of a partially ordered set $P$.

Let $P:=(\{1,\ldots, n\},\prec)$ be a partially ordered set, where $i\prec j$ implies $i<j$.
The Stanley Reisner ring of the order complex $\Delta=\Delta(P)$ is given by 
\[A:=k[\Delta]=k[x_i,\,i\in P]/\langle x_ix_j\mbox{ with $i<j$ and $i\not\prec j$}\rangle.\]

We now define a sequence of regular languages $L_i$ over the alphabet 
$\Gamma_i:=\{x_i,\ldots, x_n\}$:
\begin{enumerate}\label{reglang}
\item $x_ix_j\in L_i$ for all $i<j$ and $i\not \prec j$,\label{iter1}
\item $x_ix_{j_1}\cdots x_{j_l}\in L_i$ if $x_ix_{j_1}\cdots x_{j_{l-1}}\in L_i$ and 
  $i<j_r$ for all $r=1,\ldots l$ and either
  \begin{enumerate}
    \item $j_{l-1}\not \prec j_l$ or
    \item $x_ix_{j_1}\cdots x_{j_{l-2}}x_{j_l}\in L_i$ and $j_l<j_{l-1}$.\label{iter2}
  \end{enumerate}
\end{enumerate}
Let $f_i(x,t):=\sum_{w\in L_i}t^{|w|}\; w$ be the word counting function of $L_i$.

Corollary 3.8 and Corollary 3.9 of \cite{jolwe} 
imply the following theorem:
\begin{Thm}\label{attTupel}The Poincar\'e-Betti series of $A$ is given by:
\[\poinr{k}{A}:=\prod_{i=1}^n (1+t\, x_i)\;\;\prod_{i=1}^n(1+F_i(x,t))=
   \prod_{i=1}^n\frac{1+t\, x_i}{1-f_i(x,t)},\]
where $F_i(x,t):=\frac{f_i}{1-f_i(x,t)}$.
\end{Thm}
We only have to calculate the word counting functions $f_i$. Since the language $L_n$ is empty, it
follows that $f_n:=0$. 
We construct recursively non-deterministic finite automata $A_i$ such that the
language $L(A_i)$ accepted by $A_i$ is $L_i$
 (for the basic facts on deterministic finite automata we use here \cite{automat}). 
We assume that 
$A_j$ is defined for all $j>i$. Let $A_j^+$ be the automaton which accepts the language 
$L_j^+\cup \{w\,x_j\mbox{ with }w\in L_j^*\}$, where
\begin{eqnarray*}
L^+&:=&\big \{w_1\circ\ldots\circ w_i~\big|~i\in \N\setminus\{0\}~\mbox{and}~
   w_j\in L,~j=1,\ldots,i\big\},\\
L^*&:=&L^+\cup\{\varepsilon\}=\big \{w_1\circ\ldots\circ w_i~\big|~i\in \N~\mbox{and}~
  w_j\in L,~j=1,\ldots,i\big\},
\end{eqnarray*}
where $\circ$ denotes the concatenation and $\varepsilon$ is the empty word.
It follows that the word counting function of 
$L(A_j^+)$ is given by $\frac{f_j+t\,x_j}{1-f_j}$.\\
We now construct $A_i$:
\begin{enumerate}
\item[$\triangleright$]From the starting state we go to the state $i$ if we read the letter 
  $x_i$, otherwise we reject the input word.
\item[$\triangleright$]From the state $i$ we can switch by reading the empty word to the 
  state $j$, which represents the automaton $A_j^+$, if $i<j$ and $i\not \prec j$. 
  We then accept if $A_j^+$ accepts.
\item[$\triangleright$]Now assume we have the transitions $i\to j_1$ and $i\to j_2$ with 
  $j_1<j_2$. Because of condition (\ref{iter2}) we can switch by reading the empty 
  word from state $j_2$ to state $j_1$. 
\item[$\triangleright$] Assume that we have the transition $i\to j_2$ and we do not have the
  transition $i\to j_1$, with $j_1<j_2$. This means $i \prec j_1$ and $i \not \prec j_2$. 
  Therefore, we must have $j_1 \prec j_2$, otherwise we get a contradiction to the transitivity 
  of the order in $P$. It follows by condition (\ref{iter1}) 
  that we can switch by reading the empty word from state $j_2$ to $j_1$.
\end{enumerate}
It is clear that $A_i$ accepts the language $L_i$. Since the state $j$ represents the 
automaton $A_j^+$, we get a recursion for the word counting functions:
\begin{Lem}
For the word counting functions $f_i$ we get the following recursion:
\begin{eqnarray*}
f_n&:=&0,\\
f_i&:=&t\,x_i\;\sum_{i<j\atop i\not \prec j}\frac{f_j+t\,x_j}{1-f_j}\prod_{r=i+1}^{j-1}\frac{1+t\,x_j}{1-f_j}.
\end{eqnarray*}
\end{Lem}
\begin{proof}
The state $j$ represents the automaton $A_j^+$ with word counting function 
$\frac{f_j+t\,x_j}{1-f_j}$. 
By the argumentation above we have $j\to \nu $ for all $\nu =i+1,\ldots, j-1$ if we have $i\to j$. 
Since we accept when the automaton  $A_j^+$ accepts, we get the desired recursion.
\end{proof}
By standard facts on regular languages the functions $f_i$ are rational functions, but we 
want to have an expression of the Poincar\'e-Betti series by polynomials:
\begin{Lem}\label{rekursion f}For the rational functions $f_i$ we have:
\[f_i:=\frac{w_i}{\D 1-\sum_{r=i+1}^n w_r},\]
where $w_i$ are polynomials and $w_n=0$.
\end{Lem}
{\begin{proof} We prove it by induction: $w_n$ is a polynomial and we have 
$f_n=\frac{w_n}{1-0}$.\\
We now assume that $f_j$ satisfies the desired condition for all $j>i$. Then
{\disbreak
\begin{eqnarray*}
\lefteqn{f_i=t\,x_i\,\sum_{i<j\atop x_ix_j\in\aaf}\frac{t\,x_j+f_j}{1-f_j}\prod_{r=i+1}^{j-1}\frac{1+t\,x_r}{1-f_r}}\\
&=&t\,x_i\,\sum_{i<j\atop x_ix_j\in\aaf}\frac{t\,x_j+\frac{w_j}{1-\D\sum_{r>j}w_r}}{1-\frac{w_j}{1-\D\sum_{r>j}w_r}}
   \prod_{r=i+1}^{j-1}\frac{1+t\,x_r}{1-\frac{w_r}{1-\D\sum_{l>r}w_l}}\\
&=&t\,x_i\,\sum_{i<j\atop x_ix_j\in\aaf}\frac{t\,x_j\left(1-\D\sum_{r>j}w_r\right)+w_j}{1-\D\sum_{r\ge j}w_r}
    \left(\prod_{r=i+1}^{j-1}(1+t\,x_r)\right)\left(\prod_{r=i+1}^{j-1}\frac{1-\D\sum_{l> r}w_l}{1-\D\sum_{l\ge r}w_l}\right)\\
&=&t\,x_i\,\sum_{i<j\atop x_ix_j\in\aaf}\frac{t\,x_j\left(1-\D\sum_{r>j}w_r\right)+w_j}{1-\D\sum_{r\ge j}w_r}
    \left(\prod_{r=i+1}^{j-1}(1+t\,x_r)\right)\frac{1-\D\sum_{l> j-1}w_l}{1-\D\sum_{l\ge i+1}w_l}\\
&=&t\,x_i\,\sum_{i<j\atop x_ix_j\in\aaf}\left(w_j+t\,x_j-t\,x_j\,\D\sum_{r>j}w_r\right)
    \left(\prod_{r=i+1}^{j-1}(1+t\,x_r)\right)\frac{1}{1-\D\sum_{l\ge i+1}w_l}\\
&=&\frac{w_i}{1-\D\sum_{l\ge i+1}w_l}
\end{eqnarray*}}
with
\[w_i:=t\,x_i\,\sum_{i<j\atop x_ix_j\in\aaf}\left(w_j+t\,x_j-t\,x_j\sum_{r>j}w_r\right)
    \left(\prod_{r=i+1}^{j-1}(1+t\,x_r)\right).\]
By induction, $w_r$ is for $r>i$ a polynomial and therefore $w_i$ is a 
polynomial.
\end{proof}
\begin{Cor}The Poincar\'e-Betti series of $A$ is given by:
\[\poinr{k}{A}:=\prod_{i=1}^n (1+t\, x_i)\frac{1}{1-w_1-\ldots-w_n}\]
with
\begin{eqnarray*}
w_n&:=&0,\\
w_i&:=&t\,x_i\,\sum_{i<j\atop x_ix_j\in\aaf}\left(w_j+t\,x_j-t\,x_j\sum_{r>j}w_r\right)
    \left(\prod_{r=i+1}^{j-1}(1+t\,x_r)\right).
\end{eqnarray*}
\end{Cor}
\begin{proof}
The result is a direct consequence of Lemma \ref{rekursion f} and 
Theorem \ref{attTupel}.
\end{proof}
\medskip

We now solve the recursion of $w_i$. For this, we introduce a directed graph $G=(V,E)$ with 
vertex set $V=\{1,\ldots, n\}$ and two vertices $i,j$ are joined 
(i.e. $i\mapsto j$) if $i<j$ and $i\not \prec j$.
We write $G\big|_{i_1,\ldots, i_\nu }$ for the induced subgraph on the vertices $i_1,\ldots, i_\nu $.\\
For a sequence $1\le i_1<\ldots<i_\nu \le n$ we define
\begin{eqnarray*}
d(i_1,\ldots, i_\nu )&:=&\#\{\mbox{paths from $i_1$ to $i_\nu $ in $G\big|_{i_1,\ldots, i_\nu }$}\},\\
c(i_1,\ldots, i_\nu )&:=&\sum_{\begin{array}{c}\sst 0=a_0<a_1<\ldots<a_r=\nu \\\sst a_{i+1}-a_i\ge 2\\
    \sst r\ge 1\end{array}}(-1)^r\, d(i_{a_0+1},\ldots, i_{a_1})\cdots d(i_{a_{r-1}+1},\ldots, i_{a_r}).
\end{eqnarray*}
Note that a path counted by $d(i_1,\ldots, i_\nu )$ needs not to pass through 
all vertices $i_1,\ldots, i_\nu $.
\medskip

With this notation we get
\begin{Cor}\label{hibipoin}The Poincar\'e-Betti series of $A$ is given by:
\[\poinr{k}{A}:=\prod_{i=1}^n (1+t\, x_i)\frac{1}{W(t,\ul{x})}\]
with
\[W(t,\ul{x})=1+\sum_{1\le i_1<\ldots<i_\nu \le n\atop \nu \ge 2}c(i_1,\ldots, i_\nu )\,t^\nu \,x_{i_1}\cdots x_{i_\nu }.\]
\end{Cor}
\begin{proof}
The result follows if one solves the recursion of the $w_i$'s and collects the coefficients
of the monomials $x_{i_1}\cdots x_{i_\nu }$.
\end{proof}
In order to prove property \propP~, we give a bijection between the paths in 
$G\big|_{i_1,\ldots, i_\nu }$ and the sting-chains:
\begin{Lem}\label{bijektion}For any sequence $1\le i_1<\ldots< i_\nu \le$ there exists a 
bijection between the paths from $i_1$ to $i_\nu $ in $G\big|_{i_1,\ldots, i_\nu }$ and the 
sting-chains $I$ with $\lcm(I)=x_{i_1}\cdots x_{i_\nu }$.
\end{Lem}
\begin{proof}
We consider the path $i_1\to j_2\to j_3\to\ldots\to j_r\to i_\nu $. To this path, we associate 
the set $I:=\{x_{i_1}x_{j_2}, x_{j_2}x_{j_3},\ldots, x_{j_r}x_{i_\nu }\}$. Now we define the stings:
Assume $j_r<i_{l_0},\ldots, i_{l_1}<j_{r+1}$. 
Then we must have either $j_r\not \prec i_s$ or $i_s\not \prec j_{r+1}$ for all 
$s=l_0,\ldots, l_1$ (otherwise we would have a contradiction to $j_r\not \prec j_{r+1}$). 
This implies
\[\{x_{j_r}x_{i_s}, x_{i_s}x_{j_{r+1}}\}\cap \aaf \neq\emptyset\mbox{ for all }s=l_0,\ldots, l_1.\]
If $x_{j_r}x_{i_s}\in\{x_{j_r}x_{i_s}, x_{i_s}x_{j_{r+1}}\}\cap\aaf$, we choose $x_{j_r}x_{i_s}$,
otherwise we choose $x_{i_s}x_{j_{r+1}}$. With this choice we get that $I$ satisfies condition 
(\ref{spikecond1}) and (\ref{spikecond2}) of Definition \ref{def_sting}.
By construction we have $\lcm(I)=x_{i_1}\cdots x_{i_\nu }$.\\[+1mm]
If we start with a sting-chain $I$ with $\lcm(I)=x_{i_1}\cdots x_{i_\nu }$, then by definition 
there exist monomials
$x_{i_1}x_{j_2},x_{j_2}x_{j_3},\ldots, x_{j_r}x_{i_\nu }\in I$. This sequence defines a path
$i_1\mapsto j_2\mapsto\ldots\mapsto j_r\mapsto i_\nu $. Since both constructions are 
inverse to each other, the assertion follows.
\end{proof}
It follows:
\begin{eqnarray}
W(t,\ul{x}):=1+\sum_{I\in\Bc}(-1)^{cl(I)}m_It^{cl(I)+|I|},\label{gleichung}
\end{eqnarray}
where $\Bc$ is the set of chains of sting-chains, defined in Section \ref{sec_taylor}.

We now can prove property \propP~ and \propH~ for the ring $A=k[\Delta]$:
\begin{Thm}\label{delta}Let $P$ be a partially ordered set and $\Delta$ 
the order complex of $P$. The multigraded Poincar\'e-Betti and Hilbert series of the Stanley 
Reisner ring $A=k[\Delta]=S/\aaf$ are given by:
\begin{eqnarray*}
\poinr{k}{A}&:=&\frac{\D \prod_{i\in P} (1+t\, x_i)}{W(t,\ul{x})},\\[+2mm]
\on{Hilb}_A(\ul{x},t)&:=&\frac{W(-t,\ul{x})}{\D \prod_{i\in P} (1-t\, x_i)},
\end{eqnarray*}
where 
\begin{eqnarray*}
W(t,\ul{x})&=&1+\sum_{I\not\in\Mmc}(-1)^{cl(I)}\,m_I\,t^{cl(I)+|I|}\\
&=&1+\sum_{I\not\in\Mmc_1}(-1)^{cl(I)}\,m_I\,t^{cl(I)+|I|}\\
&=&1+\sum_{I\in\Bc}(-1)^{cl(I)}\,m_I\,t^{cl(I)+|I|}\\
&=&1+\sum_{I \on{\nbc-set}}(-1)^{cl(I)}\,m_I\,t^{cl(I)+|I|}
\end{eqnarray*}
with $\Mmc=\Mmc_1\cup\Mmc_2$ a standard matching on the Taylor resolution 
$T_\bullet$ of $\aaf$. 
\end{Thm}
\begin{proof}
The assertion is a direct consequence of Corollary \ref{sum_koszul}, 
Corollary \ref{hibipoin}, and Equation (\ref{gleichung}).
\end{proof}

%
%
\section{Proof of Conjecture \ref{poin-conj} for Several Classes of 
Algebras $A$}\label{sec_proof}
In this section we prove Conjecture \ref{poin-conj} in some special cases. 
In the first subsection, we prove the conjecture for algebras $A$ for which the
Koszul homology is an $M$-ring - a notion introduced by Fr\"oberg \cite{FR}. If in 
addition the minimal resolution of $\aaf$ has the structure of a differential-graded 
algebra, we prove property \propP~ for $A$.

In the second subsection, we prove Conjecture \ref{poin-conj} for all Koszul 
algebras. Note that this gives another proof that for a partially ordered set $P$ the 
Stanley Reisner ring $A=k[\Delta(P)]$ satisfies property \propP~ and \propH.

In the last subsection, we outline an idea for a proof of Conjecture \ref{poin-conj}
 in general.
%
%
\subsection{Proof for Algebras $A$, with $H_\bullet(K^A)$ is an M-ring}\lbreak
The first class for which we can prove Conjecture \ref{poin-conj} uses a 
theorem by Fr\"oberg \cite{FR}. We use the notation of Fr\"oberg:
\begin{Def}
A $k$-algebra $R$ isomorphic to a (non-commutative) polynomial ring $k\langle X_1,\ldots, X_r\rangle$ 
divided by an ideal $\rr$ of relations  is called 
\begin{enumerate}
\item a weak M-ring if $\rr$ is generated by relations of the following types:
  \begin{enumerate}
    \item the (graded) commutator $[X_i,X_j]=0$,
    \item $m=0$, where $m$ is a monomial in $X_i$.
  \end{enumerate}
\item an M-ring if if $\rr$ is generated by relations of the following types:
  \begin{enumerate}
    \item the (graded) commutator $[X_i,X_j]=0$,
    \item $m=0$ with $m$ a quadratic-monomial in $X_i$.
  \end{enumerate}
\end{enumerate}
\end{Def}
Now we assume that $H(K_\bullet)$ is an M-ring 
and $\Mmc$ is a standard matching. Let $R'':=k\langle Y_I, I\not\in\Mmc, 
cl(I)=1\rangle/\rr''$ be the non-commutative polynomial ring divided by an ideal $\rr''$,
where $\rr''$ is generated by the following relations:
\[Y_IY_J=(-1)^{\deg_t(Y_IY_J)}Y_JY_I,~\mbox{if}~\left\{\begin{array}{c}\gcd(m_I,m_J)=1~
    \mbox{and}~I\cup J\not\in\Mmc\\
    \mbox{for all $I,J\not\in\Mmc$ with $cl(I)=cl(J)=1$}.\end{array}\right.\]
In the notion of Fr\"oberg, $R''\otimes R'$ is the MM-ring belonging to the $M$-ring 
$R'\simeq H(K_\bullet)$.
Each literal $Y_I$ has two degrees: the total degree 
$|Y_I|:=|I|+1$ and the multidegree $\deg(Y_I):=\alpha$, with $x^\alpha = m_I$.

We define $F_\bullet:=R''\otimes_k K^A_\bullet$. Since $K^A_\bullet$ is an $A$-module, 
$F_\bullet$ is a free graded $A$-module with $\deg(m\otimes n):=\deg^{R''}_t(m)+\deg^{K^A_\bullet}_t(n)$.
Let $F_i$ be the homogeneous part of degree $i$.
The next theorem proves Conjecture \ref{poin-conj} in our situation.
\begin{Thm}\label{hauptsatz} Let $\Mmc$ be a standard matching. Assume 
$H(K_\bullet)$ an M-ring.\\
If there exists a homomorphism $s:H_\bullet(K^A)\to Z_\bullet(K^A)$, 
such that $\pi\circ s=\id_{H_\bullet(K^A)}$,
then $A$ satisfies Conjecture \ref{poin-conj}.
\end{Thm}
\begin{Cor}Under the assumptions of Theorem \ref{hauptsatz} the algebra $A$ has 
properties \propP~ and \propH.\qed
\end{Cor}
\begin{proof}[Proof of Theorem \ref{hauptsatz}]
Theorem \ref{mring} verifies the conditions for Theorem 3 in \cite{FR}.
In the proof of this theorem, Fr\"oberg shows that $F_\bullet$ defines a minimal free
resolution of $k$ as an $A$-module. By Theorem \ref{mring} the homology of the Koszul 
complex is isomorphic to the ring $R'/\rr'$. Since $H_\bullet(K^A)$ is an $M$-ring, it 
follows that the ideal $\rr'$ is generated in degree two. The construction of the 
ideal $\rr'$ implies that every standard matching ends after the second sequence.
In the second sequence of $\Mmc$, we have that $I\to J\in\Mmc_2$ satisfies 
$cl(I)=cl(J)-1$ and $|I|=|J|+1$.
Now let $I\to J\in\Mmc_2$ with $cl(I)=1$ and $cl(J)=cl(J_1)+cl(J_2)=2$. 
The difference between the ring $R''$  and the ring $R$ is that in $R$ we have 
a variable $Y_I$ and the variables $Y_{J_1},Y_{J_2}$ commute. In
the ring $R''$ the variables $Y_{J_1},Y_{J_2}$ do not commute and the variable 
$Y_I$ is omitted. Identifying $Y_{J_1}Y_{J_2}\in R''$ with 
$Y_{J_1}Y_{J_2}\in R$ and $Y_{J_2}Y_{J_1}\in R''$ with $Y_I\in R$
gives an isomorphism as $k$-vectorspaces of $R$ and $R''$.
The property $cl(I)=cl(J)-1$ and $|I|=|J|+1$ proves that this isomorphism 
preserves the degrees, and we are done.
\end{proof}
The theorem includes the theorem by Charalambous and Reeves since in their case every 
standard matching is empty and Charalambous and Reeves proved the existence of 
the map $s:H_\bullet(K^A)\to Z_\bullet(K^A)$:
\begin{Cor}[\cite{CH}] If the Taylor resolution of $\aaf$ is 
minimal, then $A=S/\aaf$ satisfies Conjecture \ref{poin-conj}.\qed
\end{Cor}
Note that $H_\bullet(K^A)\iso R'$ carries three gradings. Let $u\in R'$ with
$u=Y_{I_1}\cdots Y_{I_r}$. Then we have $\gcd(m_{I_j},m_{I_{j'}})=1$,
for $j\neq j'$, and $I_1\cup \ldots\cup I_r\not\in\Mmc$ 
(otherwise $u\in\rr'$). We set
\begin{eqnarray*}
\deg(u)&=&\alpha~\mbox{if}~\ul{x}^\alpha=m_{I_1}\cdots m_{I_r}
                                        =m_{I_1\cup \ldots\cup I_r},\\
\deg_t(u)&=&r=cl(I_1\cup \ldots\cup I_r),\\
|u|&=&|I_1|+\ldots +|I_r|=|I_1\cup \ldots\cup I_r|.\\
\end{eqnarray*}
It follows:
\[H_\bullet(K^A)\iso R'=\bigoplus_{\alpha\in\N^n\atop i,j\ge 0}R'_{\alpha,i,j}=
\bigoplus_{\begin{array}{c}I\not\in\Mmc\\\deg_t(I)=i\\|I|=j\end{array}}k~Y_I,\]
where $Y_I=Y_{I_1}\cdots Y_{I_r}$ if $cl(I)=r$ and 
$\gcd(m_{I_j},m_{I_{j'}})=1$, for $j\neq j'$.\\
Fr\"oberg proved that in the case where $H_\bullet(K^A)$ is an M-ring and
the minimal resolution of $\aaf$ has the structure of a differential-graded 
algebra we have:
\[\poinr{k}{A}=\frac{\hilb{K_\bullet\otimes_A k}}{\Hilb_{H_\bullet(K^A)}(x,-t,t)}
=\prod_{i=1}^n(1+t~x_i)\frac{1}{\Hilb_{H_\bullet(K^A)}(x,-t,t)}.\]
Therefore, we only have to calculate the Hilbert series 
$\Hilb_{H_\bullet(K^A)}(x,-t,t)$:
\begin{eqnarray*}
\Hilb_{H_\bullet(K^A)}(x,-t,t)&=&\sum_{\alpha\in\N^n\atop i,j\ge 0}
\dim_k(R'_{\alpha,i,j})~\ul{x}^\alpha~(-t)^i~t^j\\
&=&\sum_{I\not\in\Mmc}m_I~(-t)^{cl(I)}~t^{|I|}\\
&=&1~/~\Hilb_R(\ul{x},1,t).
\end{eqnarray*}
The last equation follows from Lemma \ref{secondseq}
since if $H_\bullet(K^A)$ 
is an M-ring, every standard matching ends after the second sequence. 
It follows:
\begin{Cor}
If $H_\bullet(K^A)$ is an M-ring and the minimal resolution 
of $\aaf$ has the structure of a differential-graded algebra, then
$A$ has property \propP.\qed
\end{Cor}
%
%
%
\subsection{Proof for Koszul Algebras}
In this subsection we give the proof of Conjecture \ref{poin-conj} for Koszul 
algebras $A=S/\aaf$. Note that since $\aaf$ is monomial, this is equivalent to the 
fact that $\aaf$ is generated in degree two. 
We assume in addition that $\aaf$ is squarefree. This is no restriction since
via polarization we can reduce the calculation of the Hilbert and Poincar\'e-Betti series
of $S/\aaf$ to the calculation of the series for $S/\bbf$ for a squarefree ideal 
$\bbf\subg S$.
\begin{Thm}\label{koszul}Let $A=S/\aaf$ be the quotient algebra of the 
polynomial ring and a squarefree monomial ideal $\aaf$ generated by monomials of 
degree two and $\Mmc=\Mmc_1\cup\Mmc_2$ a standard matching of $\aaf$.
Then $A$ satisfies Conjecture \ref{poin-conj}.
\end{Thm}
\begin{Cor}\label{hilb_koszul}The multigraded Poincar\'e-Betti and Hilbert 
series of Koszul algebras $A=S/\aaf$ for a squarefree monomial ideal 
$\aaf\subg S$ are given by:
\begin{eqnarray*}
\poinr{k}{A}&:=&\frac{\D \prod_{i\in P} (1+t\, x_i)}{W(t,\ul{x})},\\[+2mm]
\on{Hilb}_A(\ul{x},t)&:=&\frac{W(-t,\ul{x})}{\D \prod_{i\in P} (1-t\, x_i)},
\end{eqnarray*}
where 
\begin{eqnarray*}
W(t,\ul{x})&=&1+\sum_{I\not\in\Mmc}(-1)^{cl(I)}\,m_I\,t^{cl(I)+|I|}\\
&=&1+\sum_{I\not\in\Mmc_1}(-1)^{cl(I)}\,m_I\,t^{cl(I)+|I|}\\
&=&1+\sum_{I \on{\nbc-set}}(-1)^{cl(I)}\,m_I\,t^{cl(I)+|I|}.
\end{eqnarray*}
\end{Cor}
\begin{proof} The assertion follows directly from Theorem \ref{koszul}, the 
standard matching for ideals generated in degree two given in Section 
\ref{sec_taylor}, and the fact that, in this case, 
every standard matching ends after the second sequence.
\end{proof}
Note that if $\aaf\subg S$ is any ideal with a quadratic Gr\"obner basis, this 
corollary gives a form of the multigraded Hilbert and Poincar\'e-Betti series of 
$A=S/\aaf$ since, in this case, the series coincide with the series of 
$S/\ini(\aaf)$.
\begin{proof}[Proof of Theorem \ref{koszul}]
In this proof we sometimes consider the variables $x_1,\ldots, x_n$ as elements of the 
polynomial ring $S$ and sometimes as letters. In the second case the variables 
do not commute and we consider words over the alphabet $\Gamma:=\{x_1,\ldots,x_n\}$.
It will be clear from the context if we consider $w$ as a monomial in $S$ or as 
a word over $\Gamma$. For example, if we write $w\in\aaf$ or $x_i\mid w$, we 
see $w$ as a monomial.

For $j=1,\ldots, n$, let $\Llc_j$ be the sets of words
$x_{i_1}x_{i_2}\cdots x_{i_r}$, $r\ge 2$, over the alphabet $\{x_1,\ldots, x_n\}$,
 such that
\begin{enumerate}
\item $i_1=j<i_2,\ldots,i_r$,
\item for all $2\le l\le r$ there exists an $1\le l'<l$ such that 
  $x_{i_{l'}}x_{i_l}\in\aaf$ and $i_t>i_l$ for all $l'<t<l$.
\end{enumerate}
We define 
\[\Llc:=\left\{ w_{i_1}\cdots w_{i_r}~\left|~\begin{array}{c} i_1>\ldots>i_r\\ 
w_{i_j}\in\Llc_{i_j}, j=1,\ldots,r\end{array}\right.\right\}.\]

Note that here the variables $x_i$ are considered as letters and do not commute.
In \cite{jolwe} we construct for Koszul algebras $A$ a minimal free resolution
of $k$. The basis in homological degree $i$ in this resolution is given by the following
set (see Corollary 3.9 of \cite{jolwe}): 
\[\Bc_i=\left\{e_I~{\bf w}~\left|~\begin{array}{c}I\subset \{1,\ldots,n\}\\
{\bf w}\in\Llc\\|J|+|{\bf w}|=i\end{array}\right.\right\},\]
where $|{\bf w}|$ is the length of the word ${\bf w}$.

Thus in order to prove the theorem, we have to find a bijection between the 
words ${\bf w}\in\Llc$ of length $i$ and the monomials $u\in R$ with degree 
$|u|=i$. Remember that in our case the subsets $I\not\in\Mmc_1$ are exactly the
\nbc-sets (see Section \ref{subsec_taylor_grad2}) and 
therefore the ring $R$ has the following form: 
\[R=\frac{k\langle Y_I, I~\mbox{is an \nbc-set}~, cl(I)=1\rangle}{\langle [Y_I,Y_J]~|~
\gcd(m_I,m_J)=1\rangle}.\]
We assume that the monomials $u\in R$ are ordered, i.e. if $u=Y_{I_1}\cdots Y_{I_r}$ 
and $Y_{I_j}$ commute with $Y_{I_{j+1}}$, then $\min(I_j)>\min(I_{j+1})$.

Clearly, it is enough to construct a bijection between the sets $\Llc_j$ and 
the ordered monomials $u=Y_{I_1}\cdots Y_{I_r}$, with 
$cl(I_1\cup\ldots\cup I_r)=1$ and $j=\min(I_1)<\min(I_i)$, for $i=2,\ldots,r$.

For a word $w$ over the alphabet $\{x_1,\ldots, x_n\}$ we denote by 
$x_{f(w)}$ (resp. $x_{l(w)}$) the first (resp. the last) letter of $w$, i.e.
$w=x_{f(w)}w'$ (resp. $w=w'x_{f(w)}$).

We call a word $w$ over the alphabet $\{x_1,\ldots, x_n\}$ an \nbc-word if there
exists an index $j$ such that $w\in\Llc_j$ and each variable $x_i$, $i=1,\ldots, n$, 
appears at most once in the word $w$.

The existence of the bijection follows from the following four claims.

\medskip

\noindent {\em Claim 1:} For each $j$ and each word $w\in \Llc_j$ 
  which is not an \nbc-word there exists a unique
  subdivision  of the word $w$, \[\phi_1(w):=u_1||v_1||u_2||v_2||\ldots||u_r||v_r,\]
  such that   
  \begin{enumerate}    
    \renewcommand{\labelenumi}{(\roman{enumi})}
    \item $u_1v_1\cdots u_rv_r=w$.
    \item The subword $u_i$ is either a variable or an \nbc-word in the 
      language $\Llc_{f(u_i) }$.
    \item The words $v_i$ are either the empty word $\varepsilon$ or a descending 
      chain of variables, i.e.
      $v_i=x_{j_1}\cdots x_{j_{v_i}}$ with $j_1>\ldots>j_{v_i}$.
    \item If $v_i\neq \varepsilon$ and $u_i$ is an \nbc-word, then
      \[f(u_i)\ge f(v_i)>l(v_i)>f(u_{i+1}).\]
    \item If $v_i\neq \varepsilon$ and $u_i$ is a variable, then
      \[f(u_i)< f(v_i)>l(v_i)>f(u_{i+1}).\]
    \item If $v_i=\varepsilon$ and $u_i$ is an \nbc-word, then
      \[f(u_i)\ge f(u_{i+1}).\]
    \item If $v_i=\varepsilon$ and $u_i$ is a variable, then
      \[f(u_i)< f(u_{i+1}).\]
  \end{enumerate}

\medskip

\noindent {\em Claim 2:} There exists an injective map $\phi_2$ on the subdivisions 
  of Claim 1 such that 
  \[\phi_2\big(\phi_1(w)\big):=w_1||w_2||\ldots ||w_s\]
  and for each $w_i$, $i=1,\ldots,s$, we have the following properties: 
  
  \begin{enumerate}
    \renewcommand{\labelenumi}{(\roman{enumi})}
  \item If $w_i=x_{j_1}\cdots x_{j_t}$, then for all $1\le l\le t$ there exists an index 
    $0\le l'<l$ with $x_{j_{l'}}x_{j_l}\in\aaf$ and $j_\nu >j_l$ 
    for all $l'<\nu <l$.\label{con_1}
  \item In each word $w_i$, each variable $x_1,\ldots, x_n$ appears at 
    most once.\label{con_2}  
  \item $w_i$ is not a variable.\label{con_3}
  \item There exists an index $t$ such that $x_t\mid w_1\cdots w_{i-1}$ 
    and $x_tx_{f(w_i)}\in\aaf$ and either $x_{f(w_i)}\mid w_1\cdots w_{i-1}$ 
    or $t>f(w_i)$.\label{con_4}  
  \item For all $x_j\mid w_i$, $j<f(w_i)$, and $x_t\mid w_1\cdots w_{i-1}$ with
    $x_tx_j\in\aaf$, we have $t<j$.\label{con_5}
  \item If $\gcd(w_i,w_{i+1})=1$, then $f(w_i)>f(w_{i+1})$.\label{con_6}
  \end{enumerate}

\medskip

\noindent {\em Claim 3:} There exists an injection $\phi_3$ between 
the sequences $\phi_2\phi_1\big(\Llc_j\big)$ from Claim 2 and the 
sequences $w_1||w_2||\ldots ||w_s,$, satisfying, in addition to the 
conditions from Claim 2, the following properties:
\begin{enumerate} 
  \renewcommand{\labelenumi}{(\roman{enumi})}
  \item There exists an $j<i$ such that $\gcd(w_i,w_j)\neq 1$.\label{ord-con-1}
\end{enumerate}

\medskip

\noindent{\em Claim 4:} For each $j$ there is a bijection 
  \[\phi_4:\phi_3\phi_2\phi_1\Big(\Llc_j\Big)\to
  \left\{Y_{I_1}\cdots Y_{I_r}~\left|~\begin{array}{c}cl(I_1\cup\ldots\cup I_r)=1~\mbox{and}\\
    j=\min(I_1)<\min(I_i)~\mbox{, for $i=2,\ldots,r$}\\
    Y_{I_1}\cdots Y_{I_r}~\mbox{ordered}\end{array}\right.\right\}\]

Since $\phi_1,\ldots,\phi_3$ are injections and $\phi_4$ is a bijection, the composition 
$\phi_4\phi_3\phi_2\phi_1$ is the desired map.\\[+2mm]

\noindent{\em Proof of Claim 1.} 
Let $x_{j_1}\cdots x_{j_r}\in\Llc_j$, for 
some $j$, which is not an \nbc-word. Then we have the following uniquely defined subdivision:
\[\underbrace{x_{i_1}x_{i_2}\cdots x_{i_{j_0-1}}}_{i_2>\ldots>i_{j_0-1}}~||~
\underbrace{x_{i_{j_0}}\cdots x_{i_{j_1-1}}}_{\in\Llc_{j_0}\atop i_{j_0-1}>i_{j_0}}~||~
\underbrace{x_{i_{j_1}}\cdots x_{i_{j_2-1}}}_{i_{j_1}>\ldots>i_{j_2-1}\atop i_{j_0}\ge i_{j_1}}~||~
\underbrace{x_{i_{j_2}}\cdots x_{i_{j_3-1}}}_{\in\Llc_{j_2}\atop i_{j_2-1}>i_{j_2}}~||~
\cdots.
\]
The first part $x_{i_1}x_{i_2}\cdots x_{i_{j_0-1}}$ we split again into
\[u_1||v_1 := x_{i_1}||x_{i_2}\cdots x_{i_{j_0-1}}.\]
Thus, we get the subdivision
\[u_1~||~v_1~||~u_2~||~v_2~||\ldots||~u_{s_1}~||~v_{s_1},\]
where $u_1$ is a variable, $v_i$ are the monomials of the descending chains of 
variables (note that $v_i=\varepsilon$ is possible) and the words $u_i$, $i\ge 2$, are 
words in $\Llc_{f(u_i)}$.
If all $u_i$ are \nbc-words, we are done. But in general, it is not the case. 
Therefore, we define the following map $\varphi$:
For an \nbc-word $w$ we set $\varphi(w):=w$.
If $w$ is not an \nbc-word, we construct the above subdivision and set
\[\varphi(w):=u_1~||~v_1~||~\varphi(u_2)~||~v_2~||\ldots||~\varphi(u_{s_1})~||~v_{s_1}.\]
Since the word $w$ is of finite length the recursion, is finite and $\varphi(w)$ produces a
subdivison of the word $w$.\\
Since each $\varphi(w)$ ends with a word $v$, which is possibly the empty word 
$\varepsilon$, the $u$'s and $v$'s
do not always alternate in $\varphi(w)$.
In order to define the desired subdivision, we therefore have to modify $\varphi(w)$:
\begin{enumerate}
\item[$\triangleright$] If we have the situation $v_i||v_{i+1}$ such that $v_i, v_{i+1}$
  are descending chains of variables, possibly $\varepsilon$, then by construction we have 
  that the word $v_iv_{i+1}$ is a descending chain of variables. We replace the subdivison
  $v_i||v_{i+1}$ by the word $v_iv_{i+1}$.
\end{enumerate}
The construction implies that the resulting subdivison fulfills all desired properties.
Let $\phi_1$ be the map which associates to each word $w$ the corresponding subdivison. 
Clearly, this subdivision is unique and therefore $\phi_1$ is an injection.

\medskip

\noindent{\em Proof of Claim 2.} Let
$\phi_1(w)=u_1~||~v_1~||~u_2~||~v_2~||\ldots||~u_s~||~v_s$ be a subdivision
of Claim 1. We construct the image under $\phi_2$ by induction.
\begin{enumerate}
\item[(R)] If $f(v_s)\le f(u_s)$ and there exists a variable 
  $x_t\mid u_1v_1\cdots u_{s-1}v_{s-1}$ with $x_tx_{f(v_s)}\in\aaf$, we replace 
  $v_{s-1}$ by $v_{s-1}':=v_{s-1}x_{f(v_s)}$, else we replace $u_s$ by $u_s':=u_sx_{f(v_s)}$. 
  Finally, we replace $v_s$ by the $v_s'$ such that $v_s=x_{f(v_s)}v_s'$.
\end{enumerate}
We repeat this process until $v_s'=\varepsilon$. We get a word
\[u_1||v_1||\ldots ||u_{s-1}||v_{s-1}'||u_s',\]
such that $u_i, v_i$, for $i=1,\ldots, s-2$, and $u_{s-1}$ are as before,
$v_{s-1}'$ is a descending chain of variables and for $u_s'$ we have: 
\begin{enumerate}
\item[$(*)$]If there exist variables $x_i\mid u_s'$ with $i<f(u_s')$
  and $x_j\mid u_1v_1\cdots u_{s-1}v_{s-1}'$ such that $x_ix_j\in\aaf$,
  then $j<i$.
\end{enumerate}
Now we repeat the same process for $u_{s-1}||v_{s-1}'$.
We get a word \[u_1||v_1||\ldots ||u_{s-2}||v_{s-2}'||u'_{s-1}||u'_s,\]
such that $u_i, v_i$ are from the original decomposition and $u'_s, u'_{s-1}$
have property $(*)$.\\
We repeat this process for all words $u_i||v_i$ and we reach a sequence of words
\[\phi_{2,1}\big(\phi_1(w)\big):=u_1'||u_2'||\ldots ||u'_{s-1}||u'_s.\]
By construction this sequence satisfies the conditions 
(i), (ii), and (v).\\
Note that our construction implies that each word $u'_i$ has a unique decomposition 
$u'_i=u_i''v_i''$ such that $u_i''$ is either a variable or an \nbc-word in 
$\Llc_{f(u_i'')}$ and $v_i''$ is descending chain of variables.
Now we begin with $v_1''$ and permute the variables with respect to the rule 
(R) to the right, if necessary, and go on by induction.
It is clear that these two algorithms are inverse to each other and therefore $\phi_{2,1}$ 
is an injection onto its image.\\
In order to satisfy conditions (iii), (iv), and (vi),
we define an injective map $\phi_{2,2}$ on the image of $\phi_{2,1}$.
The composition $\phi_2:=\phi_{2,2}\phi_{2,1}$ gives then the desired map.\\
Let $\phi_{2,1}\big(\phi_1(w)\big)=u_1||u_2||\ldots ||u_{s-1}||u_s$.
Let $i$ be the smallest index such that $\gcd(u_i,u_{i+1})=1$ and 
$f(u_i)<f(u_{i+1})$. By construction the word
$u_i=u_i'v_i$ has a decomposition such that $v_i$ is a descending chain of variables
and $f(v_i)<f(u_{i+1})$ ($v_i$ was constructed by the map $\phi_{2,1}$). The word 
$u_{i+1}$ has a decomposition $u_{i+1}=u_{i+1}'v_{i+1}$ such that $u_{i+1}'$ is 
either a variable or an \nbc-word and $v_{i+1}$ a descending chain of variables.
We replace $u_i||u_{i+1}$ by the new word 
$\varphi(u_i||u_{i+1}):=u_i'u_{i+1}'c(v_iv_{i+1})$ where 
$c(v_i,v_{i+1})$ is the descending chain of variables consisting of the 
variables of $v_i$ and $v_{i+1}$.\\
We repeat this procedure until there are no words $u_i$, $u_{i+1}$ with 
$\gcd(u_i,u_{i+1})=1$ and $f(u_i)<f(u_{i+1})$.\\
It is straightforward to check that the resulting sequence 
\[\phi_{2,2}\phi_{2,1}\big(\phi_1(w)\big):=\tilde u_1||\tilde u_2||\ldots 
||\tilde u_{\tilde s-1}||u_{\tilde s}\]
satisfies all desired conditions.\\
To reverse the map $\phi_{2,2}$, we apply to each word $u_i$ the maps $\phi_1$ 
and $\phi_{2,1}$. Then it is easy to see that the sequence
\[\phi_{2,1}\phi_1(u_1)||\phi_{2,1}\phi_1(u_2)||\ldots ||\phi_{2,1}\phi_1(u_{s-1})
||\phi_{2,1}\phi_1(u_s)\]
is the preimage of $\phi_{2,2}$. Therefore, $\phi_{2,2}$ is an injection and
the map $\phi_2:=\phi_{2,2}\phi_{2,1}$ is the desired injection.

\medskip

\noindent{\em Proof of Claim 3:}
Let $\phi_2\phi_1(w)=u_1||u_2||\ldots ||u_{s-1}||u_s$ be a sequence from Claim 2.
In order to satisfy the desired condition, we construct a map $\phi_3$ similar to
$\phi_{2,2}$. Let $i$ be the largest index such that 
$\gcd(\lcm(u_1,\ldots,u_i),u_{i+1})=1$. Then it follows from Claim 2 that
$f(u_i)>f(u_{i+1})$. If we replace $u_i||u_{i+1}$ by a new word which is constructed in 
a similar way as in the map $\phi_{2,2}$, we risk to violate condition (v)
from Claim 2. Therefore, we first have to permute the word $u_{i+1}$
in the correct position. Let $l<i+1$ be the smallest index such that there 
exists an index $t>f(u_{i+1})$ with $x_t\mid u_l$ and $x_tx_{f(u_{i+1})}\in\aaf$. 
By Condition (iv) from Claim 2, such an index always exists. We replace the 
sequence $u_1||u_2||\ldots ||u_{s-1}||u_s$ by the sequence
\[u_1||\ldots||u_{l-1}||\varphi(u_l||u_{i+1})||u_{l+1}||\ldots ||u_i||u_{i+2}||\ldots||u_s,\]
where $\varphi(u_l||u_{i+1})$ is the map from the construction of $\phi_{2,2}$
of Claim 2.
Now the construction implies that all conditions of Claim 2 are still satisfied.\\
We repeat this procedure until the sequence satisfies the desired condition.\\
To reverse this procedure we reverse the map $\varphi$ with the maps $\phi_1$ and $\phi_2$
and permute the words to the right until Condition (vi) from 
Claim 2 is satisfied. It follows that $\phi_3$ is an injection onto its image.

\medskip

\noindent{\em Proof of Claim 4.} Let $\phi_3\phi_2\phi_1(w)=w_1||w_2||\ldots ||w_s$ be a
sequence from Claim 3. We now construct a bijection between these sequences
of words and the ordered monomials $Y_{I_1}\cdots Y_{I_r}$ with
$cl(I_1\cup\ldots\cup I_r)=1$ and $\min(I_1)<\min(I_j)$ for all
$j=2,\ldots, r$. We now assume:

\medskip

\noindent{\em Assumption A:}
\begin{enumerate}
\renewcommand{\labelenumi}{(\alph{enumi})}
\item For each \nbc-set $I$ and each index $i$ with $x_i\mid m_I=\lcm(I)$,
  there exists a unique word $\psi(I):=w$ such that $w=x_iw'$ and $w$ 
  satisfies conditions   (i) - (iii) from Claim 2.
\item For each word $w$ satisfying conditions (i) - (iii) 
  from Claim 2, there exists a unique \nbc-set $\varphi(w):=I$.
\end{enumerate}
In addition, the maps $\psi$ and $\varphi$ are inverse to each other.

\medskip

\noindent We now prove Claim 4:\\
Let $Y_{I_1}\cdots Y_{I_s}$ be an ordered monomial with
$cl(I_1\cup\ldots\cup I_s)=1$ and $\min(I_1)<\min(I_j)$, for $j=2,\ldots,s$.
Let $j_{I_l}$ be the smallest index $i$ such that $x_i|\lcm(I_l)$ and either
\begin{itemize}
    \item there exists a variable $x_t\mid w_1w_2\cdots w_{l-1}$ with $t>i$ 
      and $x_ix_t\in \aaf$
    \item or $x_i\mid \lcm(I_1,I_2,\ldots, I_{l-1})$. 
\end{itemize}
Such an index always exists since 
$\gcd(m_{I_1\cup I_2\cup\ldots\cup I_{l-1}},m_{I_l})\neq 1$.
By definition the variables $Y_I,Y_J$ commute if $\gcd(m_I,m_J)=1$.
It is easy to see that one can reorder the monomial $Y_{I_1}\cdots Y_{I_s}$,
such that if $\gcd(m_{I_i},m_{I_{i+1}})=1$, we have $j_{I_i}>j_{I_{i+1}}$.
We now construct a bijection between monomials $Y_{I_1}\cdots Y_{I_s}$ 
ordered in that way and the sequences of Claim 3.

Let $\phi_3\phi_2\phi_1(w)=w_1||w_2||\ldots ||w_s$ be a sequence of Claim 3 
and $I_j$ be the \nbc-sets corresponding to the words $w_j$. 
Then we associate to the sequence the following monomial
$$\phi_4(w_1||w_2||\ldots ||w_s):=Y_{I_1}\cdots Y_{I_s}.$$ 
Condition (i) from Claim 3 and Condition (vi) 
from Claim 2 imply that we get an ordered monomial.\\
On the other hand, consider an ordered monomial $Y_{I_1}\cdots Y_{I_s}$.
We associate to $Y_{I_1}$ the corresponding
\nbc-word $w_1$ whose front letter is $x_{\min(I_1)}$.\\
For $l=2,\ldots s$ let $w_l$ be the word corresponding to $I_l$ 
whose front letter is $x_{j_{I_l}}$.\\
It follows directly from the construction that the sequence 
$w_1||w_2||\ldots||w_s$ satisfies all desired conditions.\\
Conditions (iv) and (v) of Claim 2 imply
that both constructions are inverse to each other and 
therefore $\phi_4$ is a bijection.

\medskip

In order to finish our proof, we have to verify Assumption A.\\
To a word $w=x_{j_1}\cdots x_{j_s}$ satisfying Conditions (i) - (iii)
we associate a graph on the vertex set $V=[n]$.
The edges are constructed in the following way: We set $E:=\big\{\{j_1,j_2\}\big\}$.
For $j_s$ there exists an index $0\le l<s$ such that
$x_{j_l}x_{j_s}\in \aaf$. Let $P_{j_s}$ be the set of those indices.
Now let $l_2$ be the maximum of $P_{j_2}$. If $E\cup\big\{\{j_{l_2},j_2\}\big\}$
contains no broken circuit (with respect to the lexicographic order),
we set $E:=E\cup\big\{\{j_{l_2},j_2\}\big\}$.
Else we set $P_{j_2}:=P_{j_2}\setminus\{l_2\}$ and repeat the process.
It is clear that there exists at least one index in $P_{j_2}$ such that
the constructed graph contains no broken circuit.
We repeat this for $P_{j_3}, P_{j_4}, \ldots, P_{j_r}$. By construction
we obtain a graph which contains no broken circuit. Now graphs without 
broken circuits are in bijection with the \nbc-sets 
(define $I:=\{x_ix_j~|~\{i,j\}\in E\}$).\\
Given an \nbc- graph and a vertex $i$ such that there exist $j\in V$ with
$\{i,j\}\in E$, we construct a word $w$ satisfying Conditions (i) - (iii) by 
induction: Assume we can
construct to each graph of length $\nu$ and each vertex $i$ a word $w$ which 
satisfies the desired conditions.\\
Given a graph of length $\nu +1$ and a vertex $i$.
Let $P_i:=\{i<j~|~\{i,j\}\in E\}$ and
$E_1:=E\setminus \big\{\{i,j\}\in E~\big|~j\in P_i\big\}$. Then
$E\setminus E_1$ decomposes in $|P_i|+1$ connected components.
One component is the vertex $i$ and for each $j>i$ we have exactly
one component $G_j$ with $j\in G_j$. By induction we can construct words
$w_j$ corresponding to $G_j$. Now assume $P_i=\{j_1<\ldots<j_r\}$. We
set $w:=iw_{j_r}\cdots w_{j_1}$. Finally, we permute $x_t\in w_{j_l}$,
with $t<j_{l+1}$ to the right until it is in the correct position.\\
Let $w$ be a word constructed from a graph.
Assume there is $x_t\in w_j$ which was permuted to the right in the
word $w_{j'}$, $j<j'$. If there exists an index $l$ such that $x_l\in w_{j'}$,
$x_lx_t\in\aaf$, and $l>t$, then we would add an edge $\{l,t\}$.
But since $x_t\in w_j$ and the original graph was connected, this leads to 
a broken circuit for the constructed graph. Therefore, the edge for
the vertex $t$ has to be constructed with the corresponding index in $w_j$.
This proves that both constructions are inverse to each other.
\end{proof}
%
%
%
%
\subsection{Idea for a Proof in the General Case}
In this section we outline a program which we expect to yield a proof of
Conjecture \ref{poin-conj} in general.

The only way to prove the conjecture is to find a minimal $A$-free resolution 
of the field $k$, which in general is a very hard problem. With the Algebraic Discrete
Morse theory one can minimize a given free resolution, but one still needs a 
free resolution to start. The next problem is the connection to the minimized 
Taylor resolution of the ideal $\aaf$.\\
The Eagon complex is an $A$-free resolution of the field $k$ which has a 
natural connection to the Taylor resolution of the $\aaf$ since the modules 
in this complex are tensor products of 
$H_\bullet(K^A)\simeq T^\mathcal{M}\otimes_S k$.
The problem with the Eagon complex is that the differential is defined 
recursively.\\
In the first part of this section, we define a generalization of the Massey 
operations which gives us an explicit description of the differential of the 
Eagon complex. We apply Algebraic Discrete Morse theory to the Eagon complex.
The resulting Morse complex is not minimal in general, but it is minimal if
for example $H_\bullet(K^A)$ is an M-ring.
In order to prove our conjecture in general, one has to find an isomorphism 
between the minimized Eagon complex and the conjectured minimal resolution.
We can not give this isomorphism in general, but with this Morse complex we can
explain our conjecture.\\
For the general case, we think that one way to prove the conjecture is the 
following:
\begin{itemize}
\itemsep-0.6 ex
\item calculate the Eagon complex,
\item minimize it with the given acyclic matching,
\item find a degree-preserving $k$-vectorspaces-isomorphism to the ring 
  $K_\bullet\tensor_k R$.
\end{itemize}

\medskip

As before we fix one standard matching $\mathcal{M}$ on the Taylor resolution 
of $\aaf$. The set of cycles $\{\phi(I)\mid I\not\in\mathcal{M}\}$ is a system 
of representatives for the Koszul homology. With the product on the homology, we 
can define the following operation:\\
For two sets $J,I\not\in\mathcal{M}$ we define:
\[I\wedge J:=\left\{\begin{array}{lcl}0&,&\gcd(m_I,m_J)\neq 1\\
    0&,&\gcd(m_I,m_J)=1, I\cup J\in\mathcal{M} \mbox{ and }
    [\phi(I)][\phi(J)]=0\\I\cup J&,&[\phi(I)][\phi(J)]=[\phi(I\cup J)]
    \mbox{ and }I\cup J\not\in\mathcal{M}\\
    \sum_{L\not\in\mathcal{M}}a_L L&,&[\phi(I)][\phi(J)]=
    \sum_{L\not\in\mathcal{L}}a_L [\phi(L)]\mbox{ and }I\cup J\in\mathcal{M}.
  \end{array}\right.\]
Now we can define the function $(I,J)\mapsto g(I,J)\in K_\bullet^A$ such that
\[\partial(g(I,J)):=\phi(I)\phi(J)-\frac{m_Im_J}{m_{I\cup J}}\phi(I\wedge J).\]
By Proposition \ref{homo} this function is well defined.\\
We now define a function for three sets $\gamma(I_1,I_2,I_3)$ by:
\begin{eqnarray*}
\lefteqn{\gamma(I_1,I_2,I_3):=\phi(I_1)g(I_2,I_3)+(-1)^{|I_1|+1}g(I_1,I_2)
  \phi(I_3)}\\
&&+(-1)^{|I_1|+1}\frac{m_{I_1}m_{I_2}}{m_{I_1\cup I_2}}g(I_1\wedge I_2,I_3)
    -(-1)^{|I_1|+1}\frac{m_{I_2}m_{I_3}}{m_{I_2\cup I_3}}g(I_1,I_2\wedge I_3).
\end{eqnarray*}
It is straightforward to prove that $\partial(\gamma(I_1,I_2,I_3))=0$. If 
$\gamma(I_1,I_2,I_3)$ is a boundary for all sets $I_1,I_2,I_3$, we can define 
$g(I_1,I_2,I_3)$ such that $\partial(g(I_1,I_2,I_3))=\gamma(I_1,I_2,I_3)$.\\
Similar to the Massey-operation we go on by induction:\\
Assume $\gamma(I_1,\ldots, I_l)$ vanishes for all $l$-tuples $I_1,\ldots, I_l$, with 
$l\ge \nu -1$. Then there exist cycles $g(I_1,\ldots, I_l)$ such that 
$\partial(g(I_1,\ldots, I_l))=\gamma(I_1,\ldots, I_l)$. We then define:
\begin{eqnarray*}
\lefteqn{\gamma(I_1,\ldots,I_\nu ):=\phi(I_1)g(I_2,,\ldots, I_\nu )+
  (-1)^{\sum_{j=1}^{\nu -2}|I_j|+1}g(I_1,\ldots, I_{\nu -1})\phi(I_\nu )}\\
&&+\sum_{i=2}^{\nu -2}(-1)^{\sum_{j=1}^{i-1}|I_j|+1}g(I_1,\ldots, I_i)g(I_{i+1},\ldots, I_\nu )\\
&&+\sum_{i=1}^{\nu -2}(-1)^{\sum_{j=1}^i|I_j|+1}\frac{m_{I_j}m_{I_{j+1}}}{m_{I_j\cup I_{j+1}}}g(I_1,\ldots, I_{j-1}, I_j\wedge I_{j+1}, I_{j+2},\ldots, I_\nu )\\
&&-(-1)^{\sum_{j=1}^{\nu -2}|I_j|+1}\frac{m_{I_{\nu -1}}m_{I_\nu }}{m_{I_{\nu -1}\cup I_\nu }}g(I_1,\ldots ,I_{\nu -2},I_{\nu -1}\wedge I_\nu ).
\end{eqnarray*}
It is straightforward to prove that $\gamma(I_1,\ldots,I_\nu )$ is a cycle. 
Therefore, we get an induced operation on the Koszul homology.
Since the first three summands are exactly the summands of the Massey 
operations, we call $\gamma(I_1,\ldots,I_\nu )$ the $\nu $-th generalized Massey 
operations.

\medskip

>From now on we assume that all generalized Massey operations vanish. We then 
can give an explicit description of the Eagon complex:\\
We define free modules $X_i$ to be the free $A$-modules over 
$I\not\in\mathcal{M}$ with $|I|=i$. 
It is clear that we have $X_i\otimes_A k\simeq H_i(K^A)$.
The Eagon complex is defined by a sequence of complexes $Y^i$, with 
$Y^0=K_\bullet^A$ and $Y^n$ is defined by
\begin{eqnarray*}
Y_i^{n+1}&:=&Y^n_{i+1}\oplus Y_0^n\otimes X_i\mbox{,\hspace{1cm} $i>0$},\\
Y_0^{n+1}&=&Y_1^n.
\end{eqnarray*}
Let $Z_i(Y^s_\bullet)$ and $B_i(Y^s_\bullet)$ denote cycles and boundaries, 
respectively. The differentials $d^s$ on $Y^s$ are defined by induction. 
$d^0$ is the differential on the Koszul complex. Assume $d^{s-1}$ is defined. 
One has to find a map $\alpha$ that makes the diagram in Figure \ref{diagram_eins} commutative:
\begin{figure}[h]
\[
\xymatrix{
&&Y_0^s\otimes X_i\simeq Y_1^{s-1}\otimes X_i\ar[dll]_{\alpha}\ar[d]^{d^{s-1}}\\
Z_i(Y^s)\ar[rr]^{\mbox{\hspace{-1.4cm}}\pi}&&H_i(Y^s)\simeq B_0(Y^{s-1})\otimes X_i 
}\]
\caption{}\label{diagram_eins}
\end{figure}
One can then define $d^s:=(d^{s-1},\alpha)$.\\
The map $d^s$ satisfies $H_i(Y^s)=H_0(Y^s)\otimes X_i$ and 
$B_{i-1}(Y^s)=d^s(Y_1^s)=Z_i(Y^{s-1})$.
The first property allows us to continue this procedure for $s+1$ and the 
second gives us exactness of the following complex:
\[F_\bullet: \cdots Y_0^{s+1}\stackrel{d^s}{\longrightarrow}Y^s_0\stackrel{d^{s-1}}{\longrightarrow}Y_0^{s-1}\longrightarrow\cdots\longrightarrow Y_0^0\longrightarrow k.\]
Note that to make the diagram commutative, it is enough to define 
$\alpha(n\otimes f)$ for all generators
$n\otimes f$ of $Y_0^s\otimes X_i$ such that 
$\alpha(n\otimes f)=(m, d^{s-1}(n)\otimes f)$, with $m\in Y_{i+1}^{s-1}$ and 
the property that $d^{s-1}(m)+d^{s-1}(d^{s-1}(n)\otimes f)=0$.

\medskip

The $\nu$-th module of the complex $Y_\bullet^s$ is given by $Y_\nu ^s=K_j\otimes X_{i_1}\otimes \ldots \otimes X_{i_r}$ with $j+r+\sum_{j=1}^r i_j=\nu +s$. 
We fix an $R$-basis of $Y_\nu ^s$, by $e_L\otimes I_1\otimes\ldots\otimes I_r$ with 
$I_j\not\in\mathcal{M}$ and $e_L=e_{l_1}\wedge\ldots\wedge e_{l_t}$.
We are now able to define the maps $\alpha$:
Since all generalized Massey operations vanish, there exists elements 
$g(I_1,\ldots I_r)$ such that $\partial(g(I_1,\ldots I_r))=
\gamma(I_1,\ldots I_r)$
\begin{Lem}
Suppose that $d^{s-1}:Y_\bullet^{s-1}\to Y_\bullet^{s-1}$ is such that
\begin{eqnarray*}
\lefteqn{d^{s-1}(e_L\otimes I_1\otimes \ldots \otimes I_r)=\partial^K(e_L)\otimes I_1\otimes \ldots \otimes I_r }\\
&&+(-1)^{|L|}e_L\phi(I_1)\otimes I_2\otimes \ldots \otimes I_r\\
&&+(-1)^{|L|}\sum_{j=1}^{r-1}(-1)^{\sum_{i=1}^j |I_j|+1}\frac{m_{I_j}m_{I_{j+1}}}{m_{I_j\cup I_{j+1}}}
         e_L\otimes I_1\otimes\ldots \otimes I_j\wedge I_{j+1}\otimes\ldots\otimes I_r\\
&&+(-1)^{|L|}\sum_{j=1}^{r-1}(-1)^{\sum_{i=1}^j |I_j|+1}e_L~g(I_1,\ldots, I_{j+1})\otimes I_{j+2}\otimes\ldots\otimes I_r.
\end{eqnarray*}
If $n:=e_L\otimes I_1\otimes \ldots \otimes I_r\in Y_0^s$ and $J$ is a generator of $X_i$, we define $\alpha(n\otimes J)$ to be the map
that sends $n\otimes J$ to $(m,d^{s-1}(n)\otimes J)$ with
\begin{eqnarray*}
\lefteqn{m=(-1)^{|L|}(-1)^{\sum_{i=1}^r |I_j|+1}\frac{m_{I_r}m_J}{m_{I_j\cup J}}e_L\otimes I_1\otimes \ldots \otimes I_{r-1}\otimes I_r\wedge J}\\
&&+(-1)^{|L|}(-1)^{\sum_{i=1}^r |I_j|+1}e_L~g(I_1,\ldots, I_r,J).
\end{eqnarray*}
Then $\alpha$ makes the diagram in Figure \ref{diagram_eins} commutative.
\end{Lem}
\begin{proof}
We only have to check that $d^{s-1}(m)+d^{s-1}(d^{s-1}(n)\otimes f)=0$. 
This is a straightforward calculation and is left to the reader.
\end{proof}
\begin{Cor}
The map $d^s$ can be defined as follows:
\begin{eqnarray*}
\lefteqn{d^s(e_L\otimes I_1\otimes \ldots \otimes I_r)=\partial^K(e_L)\otimes I_1\otimes \ldots \otimes I_r }\\
&&+(-1)^{|L|}e_L\phi(I_1)\otimes I_2\otimes \ldots \otimes I_r\\
&&+(-1)^{|L|}\sum_{j=1}^{r-1}(-1)^{\sum_{i=1}^j |I_j|+1}\frac{m_{I_j}m_{I_{j+1}}}{m_{I_j\cup I_{j+1}}}
         e_L\otimes I_1\otimes\ldots \otimes I_j\wedge I_{j+1}\otimes\ldots\otimes I_r\\
&&+(-1)^{|L|}\sum_{j=1}^{r-1}(-1)^{\sum_{i=1}^j |I_j|+1}e_L~g(I_1,\ldots, I_{j+1})\otimes I_{j+2}\otimes\ldots\otimes I_r.
\end{eqnarray*}
\end{Cor}

With this corollary we get an explicit description of the Eagon resolution of 
$k$ over $A$.\\

In order to define the acyclic matching, we first use Theorem \ref{mring} to
define the Eagon complex with the ring 
$H_\bullet(K^A)\iso R'=k[Y_I\mid cl(I)=1, I\not\in\mathcal{M}]/\rr'$ instead 
of $H_\bullet$. The operation $I\wedge J$ then is nothing but the 
multiplication $Y_IY_J$ in $R'$. We write $y_I$ for the class of $Y_I$ in $R'$.

\medskip

It is clear that this complex is not minimal in general. The idea now is to 
minimize this complex via Algebraic Discrete Morse theory. It is easy to see, 
that the only invertible coefficient occurs by mapping 
$\ldots\otimes y_I\otimes y_J\otimes\ldots$ to the element
$\ldots\otimes y_Iy_J\otimes\ldots$, with $\gcd(m_I,m_J)=1$.
The idea is to match all such basis elements, with $I\wedge J=I\cup J$ and 
$I\cup J\not\in\mathcal{M}$. In order to do this, we have to define an order 
on the variables $y_I$ with $I\not\in\mathcal{M}$:
We order the sets $I$ by cardinality and if two sets have the same cardinality 
by the lexicographic order on the multidegrees $m_I, m_J$.
The monomials in $R'$ are ordered by the degree-lexicographic order.
The acyclic matching is similar to the Morse matching on the normalized Bar 
resolution (see \cite{jolwe}). Since $\mathcal{M}$ is a standard 
matching on the Taylor resolution, we know that if 
$I_1\cup I_2\cup\ldots\cup I_r\not\in\mathcal{M}$ with $cl(I_j)=1$ and 
$\gcd(m_{I_j},m_{I_{j'}})=1$ for all $j\neq j'$, then it follows that 
$I_2\cup\ldots\cup I_r\not\in\mathcal{M}$. Therefore, the following matching 
is well defined:
\[e_L\otimes y_{I_1}\otimes y_{I_2}\cdots y_{I_r}\otimes\ldots \mapsto 
e_L\otimes y_{I_1}y_{I_2}\cdots y_{I_r}\otimes\ldots,\]
where $I_1<I_2<\ldots <I_r$ and $I_1\cup I_2\cup\ldots\cup I_r\not\in\mathcal{M}$
and $cl(I_j)=1$ and $\gcd(m_{I_j},m_{I_{j'}})=1$ for all $j\neq j'$.
On the remaining basis elements we do the same matching on the second 
coordinate, and so on. The exact definition of the acyclic matching and the 
proof is given in Definition 3.1 of \cite{jolwe}.

We describe the remaining basis elements, as in \cite{jolwe}, by induction.\\
$[y_I|u_1]$ with $u_1=y_{J_1}\cdots y_{J_r}$ is called fully attached 
(see Definition 3.3 of \cite{jolwe}) if one of the 
following conditions is satisfied:
\begin{enumerate}
\item $r=1$ and $\gcd(m_I,m_{J_1})\neq 1$ or $y_I>y_{J_1}$,
\item $\gcd(m_I,m_{J_i})=1$ for all $i$ and 
  $I\cup J_1\cup\ldots\cup J_r\in\mathcal{M}$, 
  and for all $1\le i\le r$ we have 
  $I\cup J_1\cup\ldots\cup\widehat{J_i}\cup\ldots\cup J_r\not\in\mathcal{M}$.
\end{enumerate}
A tuple $[y_J|u_1|\ldots |u_r]$ is called fully attached if
$[y_J|u_1|\ldots |u_{r-1}]$ is fully attached, one of the following 
properties is satisfied and $u_r$ is minimal in the sense that there is no 
proper divisor $v_r\mid u_r$ satisfying one of the conditions below:
\begin{enumerate}
\item $u_r$ is a variable and $\gcd(m_{u_{r-1}},m_{u_r})\neq 1$,
\item $u_r, u_{r-1}$ are both variables and $u_{r-1}>u_r$,
\item $[y_J|u_1|\ldots |u_{r-2}|u_r]$ is a fully attached tuple and 
  $u_{r-1}>u_r$,
\item $u_{r-1}=y_{I_1}\cdots y_{I_t}$, $u_r=y_{J_1}\cdots y_{J_s}$ such that
  $\gcd(m_{u_{r-1}},m_{u_r})=1$ and 
  $I_1\cup\ldots\cup I_t\cup J_1\cup\ldots\cup J_s\in\Mmc$.
\end{enumerate}
Here $m_u:=\lcm(I_1\cup\ldots\cup I_r)$ if $u=y_{I_1}\cdots y_{I_r}$.

The basis of the Morse complex is given by elements $e_L|{\bf w}$, where 
${\bf w}$ is a fully attached tuple. If $H_\bullet(K^A)$ is an M-ring, 
the Morse complex is minimal since in this case the fully attached tuple 
has the form $[y_{I_1}|y_{I_2}|\cdots|y_{I_r}]$.
In order to prove Conjecture \ref{poin-conj} one has to find an isomorphism
between the fully attached tuples and the monomials in $R$.

We can not give this isomorphism in general, but we think that this Morse 
complex helps for the understanding of our conjecture:

Let $[y_{I_1}|y_{I_2}|\ldots |y_{I_r}]$ be a fully attached tuple, with
$y_{I_1}>\ldots>y_{I_r}$. We map such a tuple to the monomial
$Y_{I_1}\cdots Y_{I_r}\in R$. Clearly, this map preserves the degree.
We get a problem if $[y_J|u_1|\ldots |u_r]$ is a fully attached tuple and
$u_1=I_1\cup \ldots \cup I_r$ with $r>1$.
For example, assume $J\mapsto I_1\cup \ldots \cup I_r\in\Mmc_r$, with 
$cl(J)=cl(I_1)=\ldots=cl(I_r)=1$ and $\gcd(m_{I_j},m_{I_{j'}})=1$
for $j\neq j'$, is matched.
Assume further $y_{I_1}<\ldots<y_{I_r}$. Then 
$[y_{I_1}|y_{I_2}\cdots y_{I_r}]$ is a fully attached tuple.
We cannot map $[y_{I_1}|y_{I_2}\cdots y_{I_r}]$ to 
$Y_{I_1}Y_{I_2}\cdots Y_{I_r}$, since in $R$ the variables commute, i.e.
$Y_{I_1}Y_{I_2}\cdots Y_{I_r}=Y_{I_r}Y_{I_{r-1}}\cdots Y_{I_1}$  and 
the tuple $[y_{I_r}|y_{I_{r-1}}|\ldots |y_{I_1}]$ maps already to this element.
But we can define
\[[y_{I_1}|y_{I_2}\cdots y_{I_r}]\mapsto Y_J\in R.\]
The degree of $Y_J\in R$ is $|J|+1$ and the homological degree of 
$[y_{I_1}|y_{I_2}\cdots y_{I_r}]$ is
$$|I_1|+1+(|I_2|+\ldots +|I_r|)+1=(|I_1|+\ldots+|I_r|+1)+1= |J|+1,$$
therefore this map preserves the degree.

These facts demonstrate that the variables $Y_I$, for which $I\in\Mmc$, $cl(I)=1$, and 
$I\not\in\Mmc_1$, are necessary. We consider this as a justification of our
conjecture.

%
%
%
\section{Applications to the Golod Property of Monomial Rings}\label{sec_golod}
In this section we give some applications to the Golod property. 
Remember that a ring $A$ is 
Golod if and only if one of the following 
conditions is satisfied (see \cite{Gull}):
\begin{eqnarray}
\poinr{k}{A}=\frac{\D\prod_{i=1}^n(1+x_i\;t)}
  {\D 1-t\;\sum_{\alpha\in\mathbb{N}^n,i\ge 0}\dim_k(Tor^S_i(A,k)_\alpha) 
    x^\alpha\;t^i}.\label{golodcon}\\
\mbox{All Massey operations on the Koszul homology vanish.}
\end{eqnarray}
If an algebra satisfies property \propP, then we get in the monomial case the 
following equivalence:
\begin{Thm}\label{newequiv}If $A=S/\aaf$ satisfies property \propP~, then $A$ 
  is Golod if and only if one of the following conditions is satisfied: 
  \begin{enumerate}
    \item For all subsets $I\subset \SM(\aaf)$ with $cl(I)\ge 2$ we have 
      $I\in\Mmc$ for any standard matching $\Mmc$.
    \item The product (i.e. the first Massey operation) on the Koszul homology 
      is trivial.
  \end{enumerate}
\end{Thm}
\begin{proof}
Property \propP~ implies the equivalence of (\ref{golodcon}) and the 
first condition. Theorem \ref{mring} implies the equivalence of the first and 
the second condition. 
\end{proof}
\begin{Cor}
If $A=S/\aaf$ satisfies one of the following conditions, then $A$ is Golod 
if and only if the first Massey operation vanishes.
\begin{enumerate}
\item $\aaf$ is generated in degree two, 
\item $H_\bullet(K^A)$ is an $M$-ring and either there is a homomorphism 
  $s:H_\bullet(K^A)\to Z_\bullet(K^A)$ such that $\pi\circ s=\id_{H_\bullet(K^A)}$
  or the minimal resolution of $\aaf$ has the structure of a differential 
  graded algebra.
\end{enumerate}
\end{Cor}
\begin{proof}
In the previous section we proved property \propP~ in these cases, 
therefore the result follows from the theorem above.
\end{proof}

Recently, Charalambous proved in \cite{chara2} a criterion for generic ideals to be Golod.
Remember that a monomial ideal $\aaf$ is generic if the multidegree of two minimal
monomial generators of $\aaf$ are equal for some variable, then there is a third monomial 
generator of $\aaf$ whose multidegree is strictly smaller than the multidegree of the 
least common multiple of the other two. It is known that for generic ideals $\aaf$ the 
Scarf resolution is minimal. Charalambous proved the following proposition:
\begin{Prop}[\cite{chara2}] Let $\aaf\subg S$ be a generic ideal. $A=S/\aaf$ is Golod if 
and only if $m_Im_J\neq m_{I\cup J}$ whenever $I\cup J\in \Delta_S$ for $I,J\subset \SM(\aaf)$.\\
Here $\Delta_S$ denotes the Scarf resolution.
\end{Prop}

Assuming property \propP, our Theorem \ref{newequiv} gives a second proof of this fact:
\begin{proof}
It is easy to see that the condition
\[m_Im_J\neq m_{I\cup J}~\mbox{whenever}~I\cup J\in \Delta_S\]
is equivalent to fact that the product on the Koszul homology is trivial.
Thus, Theorem \ref{newequiv} implies the assertion.
\end{proof}

\medskip

We have the following criterion:

\begin{Lem}Let $A=S/\aaf$ with $\aaf=\langle m_1,\ldots, m_l\rangle$.
\begin{enumerate}
\item If $\gcd(m_i,m_j)\neq 1$ for all $i\neq j$, then $A$ is Golod 
  (see \cite{CH}, \cite{wel}).
\item If $A=S/\aaf$ is Golod, then $\aaf$ satisfies the $\gcd$-condition.
\end{enumerate}
\end{Lem}
\begin{proof}
If a ring $A$ is Golod, then the product on $H_\bullet(K^A)$ is trivial. This 
implies $Y_IY_J=0$ if $\gcd(m_I,m_J)= 1$. With Theorem \ref{mring} it follows
that all sets $I\cup J$ with $\gcd(m_I,m_J)=1$ are matched. In particular, all 
sets $\{m_i,m_j\}$ with $\gcd(m_i,m_j)=1$. Such a set can only be matched with 
a set $\{m_{i_1},m_{i_1},m_{i_1}\}$ with the same $\lcm$. 
But this implies that there must exist a third generator $m_r$ with 
$m_r|m_i\,m_j$.
\end{proof}
The following counterexample shows that the converse of the second statement 
is false: Let $\aaf:=\langle xy, yz, zw, wt, xt\rangle$ be the Stanley Reisner 
ideal of the triangulation of the 5-gon. It is easy to see that $\aaf$ 
satisfies the $\gcd$-condition. But $\aaf$ is Gorenstein and therefore not 
Golod. But we have:
\begin{Thm}\label{strong-gcd-cond-golod}
If $A=S/\aaf$ has property \propP~ and $\aaf$ satisfies the strong 
$\gcd$-condition, then $A$ is Golod.
\end{Thm}
\begin{proof}
We prove that $H_\bullet(K^A)$ is an $M$-ring and isomorphic as an algebra to 
the ring
\[R:=k(Y_I~\mid~I\not\in\Mmc, cl(I)=1)/\langle Y_IY_J~
\mbox{for all $I,J\not\in\Mmc_0\cup\Mmc$}\rangle,\]
where $\Mmc_0$ is the sequence of matchings constructed in Proposition 
\ref{gcd-complex} in order to obtain the 
complex $T_{\gcd}$ and $\Mmc$ is a standard matching on the complex 
$T_{\bf \gcd}$.
It follows that the first Massey operation is trivial and then Theorem 
\ref{newequiv} implies the assertion.\\
The idea is to make the same process as in Section \ref{sec_koszul} with the 
complex $T_{\bf \gcd}$ from Proposition \ref{gcd-complex} from Section 
\ref{sec_taylor} instead of the Taylor resolution $T_\bullet$. 
Since all sets $I$ in $T_{\bf \gcd}$ satisfy $cl(I)=1$, the result follows 
directly from property $\propP$.\\
Note that $\Mmc_0$ satisfies all conditions required in the proof of 
Proposition \ref{homo} except the following:
Assume $I\cup J\in\Mmc_0$ with $\gcd(m_I,m_J)=1$ and $I,J\not\in\Mmc_0$. 
Then there exists a set $\hat I$ such that $\hat I\to I\cup J\in\Mmc_0$. It follows
\[0=\partial^2(\hat I)=\partial(I\cup J)+\sum_{L\not\in\Mmc_0}a_L\;L\]
and therefore as in the proof of Proposition \ref{homo}
\[\phi(I\cup J)=\sum_{L\not\in\Mmc_0}a_L\;\phi(L)\mbox{~\hspace{5mm}for some $a_L\in k$}.\]
In the case of Proposition \ref{homo} we could guarantee 
that $cl(L)\ge cl(I\cup J)$. We can not deduce this fact here, but this is the only 
difference between $\Mmc_0\cup \Mmc$ and a standard matching on 
the Taylor resolution.
Since all sets $L$ with $cl(L)\ge 2$ are matched, we only could have
\[\phi(I\cup J)=\sum_{L\not\in\Mmc_0\atop cl(L)=1}a_L\;\phi(L)
  \mbox{~\hspace{5mm}for some $a_L\in k$}.\]
We prove that this cannot happen.
If $I\cup J$ is matched, then there exists a monomial $m$ with 
$I\cup J\cup\{m\}\to I\cup J\in\Mmc_0$.
But then, since $cl(I\cup J\setminus\{n\})\ge cl(I\cup J)\ge 2$, by the 
definition of $\Mmc_0$ any image 
$I\cup J\cup\{m\}\setminus\{n\}$ is also matched:
\[I\cup J\cup\{m\}\setminus\{n\}\to I\cup J\setminus\{n\}\in\Mmc_0.\]
This proves that the situation above is not possible and we are done.
\end{proof}

\begin{Cor}
Suppose that $A=S/\aaf$ has property \propP.~ Then $A$ is Golod if
\begin{enumerate}
\item $\aaf$ is shellable (for the definition see \cite{ekki}),
\item $\SM(\aaf)$ is a monomial ordered family (for the definition see \cite{monom}),
\item $\aaf$ is stable and $\#\supp(m)\ge 2$ for all $m\in \SM(\aaf)$,
\item $\aaf$ is $p$-Borel fixed and $\#\supp(m)\ge 2$ for all $m\in \SM(\aaf)$.
\end{enumerate}
Here $\supp(m):=\big\{1\le i\le n~\big|~x_i~\mbox{divides}~m\big\}$.
\end{Cor}
\begin{proof}
We order $\SM(\aaf)$ with the lexicographic order. Then it follows directly 
from the definitions of the ideals that $\aaf$ satisfies the strong 
$\gcd$-condition. The assertion follows then from Theorem 
\ref{strong-gcd-cond-golod}.
\end{proof}

Theorem \ref{strong-gcd-cond-golod} and the preceding Lemma give rise to the 
following conjecture:
\begin{Conj}Let $\aaf=\langle m_1,\ldots, m_l\rangle\subset S$ be a monomial 
ideal and $A=S/\aaf$.\\
Then $A$ is Golod if and only if $\aaf$ satisfies the strong $\gcd$-condition.\\
In particular: Golodness is independent of the characteristic of $k$.
\end{Conj}

\medskip

It is known that if $\aaf$ is componentwise linear, then $A$ is 
Golod (see \cite{wel}). 
One can generalize this result to the following:
\begin{Cor}
Let $\aaf$ be generated by monomials with degree $l$. 
\begin{enumerate}
\item If $\dim_k\big(Tor_i^S(S/\aaf,k)_{i+j}\big)=0$ for all $j\ge 2(l-1)$, then $A=S/\aaf$ is 
  Golod,
\item if $A$ is Golod, then $\dim_k\big(Tor_i^S(S/\aaf,k)_{i+j}\big)=0$ for all $j\ge i(l-2)+2$.
\end{enumerate}
In particular:
If $A$ is Koszul, then $A$ is Golod if and only if the minimal free resolution of $\aaf$ 
is linear.
\end{Cor}
\begin{proof} Let $I\subset \{m_1,\ldots, m_l\}$ with $cl(I)=1$ and 
  $\lcm(I)\neq\lcm(I\setminus\{m\})$ for all $m\in I$.
  Then $l+|I|-1\le \deg(I)\le (l-1)|I|+1$. 
  Now assume that $L=I\cup J\not\in\Mmc$ with $\gcd(m_I,m_J)=1$, then 
  $\deg(L)\ge 2l-2+|I\cup J|$, which is a contradiction to 
  $\dim_k\big(Tor_i^S(S/\aaf,k)_{i+j}\big)=0$ for all $j\ge 2l-2$. Therefore, the product on the 
  Koszul homology is trivial. By the same multidegree reasons 
  it follows that all Massey operations have to vanish, hence $A$ is Golod.\\
  If $A$ is Golod, then the product on $H_\bullet(K^A)$ is trivial, hence 
  (by theorem \ref{mring}) $I\not\in\Mmc$ implies $cl(I)=1$.
  But for those subsets we have $l+|I|-1\le \deg(I)\le (l-1)|I|+1$. Therefore, 
  it follows that $\dim_k\big(Tor_i^S(S/\aaf,k)_{i+j}\big)=0$ for all $j\ge i(l-2)+2$.
\end{proof}

\section{Acknowledgement}
The author would like to thank J\"urgen Herzog and Volkmar Welker 
for numerous discussions.

%
%
%

%
%
\end{document}